\documentclass[12pt,leqno]{article}
\usepackage{amsmath, amsthm, amsfonts, amssymb, color}
\usepackage{mathrsfs}
\usepackage[notcite,notref]{showkeys}

\theoremstyle{definition}

\newcommand{\scr}[1]{\mathscr #1}
\definecolor{wco}{rgb}{0.5,0.2,0.3}

\numberwithin{equation}{section} \theoremstyle{remark}

\newcommand{\ua}{\uparrow}

\title{
{\bf Gradient Estimates and Applications for SDEs in Hilbert Space  with Multiplicative Noise and Dini Continuous  Drift}\footnote{supported in part by
NNSFC(11131003) and Laboratory of  Mathematics and  Complex Systems.}}
\author{
{\bf Feng-Yu Wang \footnote{ wangfy@bnu.edu.cn,
F.-Y.Wang@swansea.ac.uk.  }}\\
 \footnotesize{School of Mathematical Sciences,
Beijing Normal
University, Beijing 100875, China}\\
 \footnotesize{Department of Mathematics,
Swansea University, Singleton Park, SA2 8PP, United Kingdom}
 }

\begin{document}
\def\tttext#1{{\normalfont\ttfamily#1}}
\def\R{\mathbb R}  \def\ff{\frac} \def\ss{\sqrt} \def\B{\mathbf B}
\def\N{\mathbb N} \def\kk{\kappa} \def\m{{\bf m}}
\def\dd{\delta} \def\DD{\Delta} \def\vv{\varepsilon} \def\rr{\rho}
\def\<{\langle} \def\>{\rangle} \def\ggm{\Gamma} \def\ggm{\gamma}
  \def\nn{\nabla} \def\pp{\partial} \def\EE{\scr E}
\def\d{\text{\rm{d}}} \def\bb{\beta} \def\aa{\alpha} \def\D{\scr D}
  \def\si{\sigma} \def\ess{\text{\rm{ess}}}
\def\beg{\begin} \def\beq{\begin{equation}}  \def\F{\scr F}
\def\Ric{\text{\rm{Ric}}} \def\Hess{\text{\rm{Hess}}}
\def\e{\text{\rm{e}}} \def\ua{\underline a} \def\OO{\Omega}  \def\oo{\omega}
 \def\tt{\tilde} \def\Ric{\text{\rm{Ric}}}
\def\cut{\text{\rm{cut}}} \def\P{\mathbb P}
\def\C{\scr C}     \def\E{\mathbb E}
\def\Z{\mathbb Z} \def\II{\mathbb I}
  \def\Q{\mathbb Q}  \def\LL{\Lambda}
  \def\B{\scr B}    \def\ll{\lambda}
\def\vp{\varphi}\def\H{\mathbb H}\def\ee{\mathbf e}

\maketitle
\begin{abstract} Consider the stochastic evolution equation in a separable Hilbert space $\H$ with a nice multiplicative noise and a locally Dini continuous drift.  We prove that   for any   initial data  the equation has a unique (possibly explosive) mild solution.   Under a reasonable condition ensuring the non-explosion of the solution, the strong Feller property of the associated Markov semigroup is  proved.  Gradient estimates and log-Harnack inequalities are derived for the associated semigroup under certain global conditions, which are new even in finite-dimensions.
\end{abstract} \noindent


\section{Introduction}

Let $(\H,\<\cdot,\cdot\>,|\cdot|)$ and $(\bar \H,\<\cdot,\cdot\>_{\bar \H},|\cdot|_{\bar \H})$ be two separable Hilbert spaces. Let $W=(W_t)_{t\ge 0}$ be a cylindrical Brownian motion on $\bar\H$ with respect to  a complete filtration probability space $(\OO,(\F_t)_{t\ge 0},\F,\P)$. More precisely, $W_t=\sum_{n=1}^\infty B_t^n \bar e_n$ for a sequence of  independent one-dimensional Brownian motions $\{B^n_\cdot\}_{n\ge 1}$ with respect to  $(\OO,(\F_t)_{t\ge 0},\F,\P)$, where $\{\bar e_n\}_{n\ge 1}$ is an orthonormal basis on $\bar \H$. Consider the following semi-linear stochastic partial differential equation on $\H$:
\beq\label{E1} \d X_t= \big\{A X_t+ B_t(X_t)+b_t(X_t)\big\}\d t + Q_t (X_t) \d W_t,\end{equation}
where $(A,\D(A))$ is a negative definite self-adjoint operator on $\H$, $B,b: [0,\infty)\times \H\to \H$ are measurable and locally bounded (i.e. bounded on bounded sets), and $Q: [0,\infty)\times\H\to \scr L (\bar\H;\H)$ is measurable, where $\scr L(\bar\H;\H)$ is the space of bounded linear operators from $\bar\H$ to $\H$.   Here, $B$ and $b$ stand  for the regular part and the singular part of the drift  respectively.

Let $\|\cdot\|$ and $\|\cdot\|_{HS}$ denote the operator norm and the Hilbert-Schmidt norm respectively, and let $\scr L_{HS}(\bar\H;\H)$ be the space of all Hilber-Schmidt operators from $\bar \H$ to $\H$. Throughout the paper, we let $A,B$ and $Q$ satisfy the following two assumptions.

   \beg{enumerate} \item[{\bf (a1)}] $(A, \D(A))$ is a negative definite self-adjoint operator on $\H$   such that $(-A)^{\vv-1}$ is of trace class for some $\vv\in (0,1);$ i.e. $\sum_{n=1}^\infty \ff 1 {\ll_n^{1-\vv}}<\infty$ for    $0<\ll_1\le\ll_2\le\cdots  $ being all eigenvalues of $-A$ counting   multiplicities.
\item[{\bf (a2)}] $B\in C([0,\infty)\times\H;\H), Q\in C([0,\infty)\times\H; \scr L(\bar\H;\H))$ such that for every $(t,x)\in [0,\infty)\times\H$,  $B_t:\H\to \H$ is local  Lipschitz continuous, $Q_t\in C^2(\H;\scr L(\bar\H;\H))$,    $(Q_tQ_t^*)(x)$ is invertible and a.e. right-continuous in $t\ge 0$, and
$$\|\nn B_t(x)\|+ \|\nn Q_t (x)\|  + \|\nn^2  Q_t(x)\|    +\|Q_t(x)\| + \|(Q_tQ_t^*)(x)^{-1}\|$$ is locally bounded
in $(t,x)\in[0,\infty)\times\H$, where $\|\nn B_t(x)\|$ stands for the local Lipschitz constant of $B_t$ at point $x$.
 \end{enumerate}

Under {\bf (a1)} and {\bf (a2)}, we  first search for minimal conditions on $b$   ensuring   the existence and pathwise uniqueness of    mild solutions to \eqref{E1}, then study gradient estimates and Harnack inequalities of the associated semigroup.

Before moving on, we briefly recall some recent   progresses  made in this direction for $\H=\bar\H$, constant  $Q$ and $B=0$.
By  {\bf (a1)} and {\bf (a2)},     the  Ornstein-Ulenbeck   process
$$ Z_t^x:= \e^{At}x +\int_0^t \e^{(t-s)A} Q\d W_s,\ \ t\ge 0, x\in \H$$  is a continuous Markov process on $\H$ having a  unqiue invariant probability measure $\mu$, see e.g. \cite{DZ}. When $Q=I$ (the identity operator), $ B=0$ and $b_t=b$ is  independent of  $t$   satisfying  a reasonable growth condition,   the existence and uniqueness of mild solutions to \eqref{E1} are proved in \cite{DR2} for $\mu$-a.e. starting points (see also \cite{RN}  for the case with an additional   gradient term). This improves  the corresponding result derived in \cite{DR1} where $b$ is   bounded. As for mild solutions to \eqref{E1} with arbitrary initial points, the existence and uniqueness have been proved in \cite{DF} when   $b$ is   H\"older continuous.

We would also like to mention that for  SDEs on $\R^d$ with a nice non-degenerate multiplicative noise, the existence and uniqueness of solutions have been proved in \cite{Zh} if the drift is in $L^{2(d+1)}_{loc}([0,\infty)\times \R^d)$. When the noise is  non-degenerate and  additive,  this condition is weakened in \cite{KR} as that the drift belongs to $L_{loc}^q([0,\infty)\to L_{loc}^p(\R^d))$ for some $p,q\in [1,\infty]$ satisfying $\ff d p+ \ff 2 q<1.$ The main idea  used in these two papers goes back to the arguments developed in \cite{V,Z} using Sobolev regularities  of the corresponding Kolmogorov equations. As already explained in e.g. \cite{DR1} that such  regularities are  not  available in infinite dimensions. Indeed, \cite{DR1,DR2} are attempts  to extend these results to infinite-dimensions by using (local) boundedness conditions to replace the local  integrability conditions.

By refining the argument developed  from \cite{DF} for additive noise, and by carefully treating the operator-valued map $Q$, 
we find that the existence and uniqueness of mild solutions can be ensured by {\bf (a1)} and {\bf (a2)} provided $Q$ is asymptotically cylindrical and $b$ is
locally Dini continuous. precisely, 
 for any $n\ge 1$, let $\pi_n:\H\to \H_n:=\text{span}\{e_1,\cdots, e_n\}$ be the orthogonal projection, where $\{e_n\}_{n\ge 1}$ is the eigenbasis of $-A$ on $\H$ corresponding to the eigenvalues $\{\ll_n\}_{n\ge 1}.$  Moreover, let
\beg{equation*}\beg{split} 
 &\D=\bigg\{\phi:  [0,\infty)\to [0,\infty)  \text{\ is\ increasing},  \phi^2 \text{\ is\ concave,} \int_0^1\ff{\phi(s)}s\d s<\infty\bigg\}.\end{split}\end{equation*}  
 We shall need   the following   condition.

\beg{enumerate}
\item[{\bf (a3)}]  $b: [0,\infty)\times\H\to\H$ is measurable and locally bounded, and for any $n\ge 1,$ there exists $\phi_n\in\D$ such that  \beq\label{BQ2}|b_t(x)-b_t(y)| \le \phi_n(|x-y|),\ \ \  t\in [0,n],\, x,y\in\H,   |x|\lor |y|\le n. \end{equation}
     Moreover, for any $x\in \H$ and $s\ge 0$,
\beq\label{QQ}\lim_{n\to \infty} \big\|Q_s(x) -Q_s(\pi_n x)\big\|_{HS}^2:=\lim_{n\to \infty} \sum_{k\ge 1} \big|\{Q_s(x) -Q_s(\pi_n x)\}\bar e_k\big|^2=0,\end{equation} where $\{\bar e_k\}$ is an orthonormal basis on $\bar\H$.
\end{enumerate}

We remark that the condition $\int_0^1\ff{\phi(s)}s\d s<\infty$ is well known as Dini condition, due to the notion of Dini continuity.
So, \eqref{BQ2} implies that $b_t$ is Dini continuous on bounded sets in $\H$, locally uniformly in $t\ge 0$. Obviously, the class   $\D$ contains $\phi(s):= \ff{K}{\log^{1+\dd}(c+s^{-1})}$ for constants $K,\dd>0$ and large enough $c\ge \e$ such that $\phi^2$ is concave.

Next, a map  $Q$ defined on $\H$ is called cylindrical if $Q(x)=Q(\pi_n x)$ holds for some $n\ge 1$ and all $x\in \H$. So, \eqref{QQ} means that $Q_s$ is asymptotically cylindrical under the Hilber-Schmidt norm, uniformly in $s\ge 0$.
We stress that   assumptions {\bf (a2)}  and {\bf (a3)} are satisfied by some infinite-dimensional models. For instance, when $\H=\bar \H$ and $Q_s(x)= Q_0+\vv \tt Q(x)$, where $Q_0\in \scr L(\bar \H;\H)$ such that $Q_0Q_0^*$ is invertible,   $\tt Q\in C_b^2(\H; \scr L(\bar \H;\H))\cap C_b (\H; \scr L_{HS}(\bar \H;\H))$  and $\vv\in\R$,  all conditions on $Q$ included in these two assumptions hold provided $|\vv|$ is small enough.

\

In general,     the mild solution (if exists) can be explosive. So, we consider   mild solutions with life times.

\beg{defn} A  continuous adapted process $(X_t)_{t\in [0,\zeta)}$ is called a mild solution to \eqref{E1} with life time $\zeta$, if $\zeta>0$ is a stopping time such that $\P$-a.s. $\limsup_{t\uparrow \zeta} |X_t|=\infty$ holds on $\{\zeta<\infty\}$ and, $\P$-a.s.
$$X_t= \e^{At}X_0 +\int_0^t \e^{(t-s)A}(B_s+b_s)(X_s)\d s +\int_0^t \e^{(t-s)A}Q_s(X_s)\d W_s,\ \ t\in [0,\zeta).$$
\end{defn}

If for any $x\in\H$, the equation \eqref{E1} has a unique mild solution $X_t^x$ with $X_0=x$ and infinite life time (i.e. the solution is non-explosive), then the associated Markov semigroup $P_t$ is defined as follows.
$$P_t f(x):= \E f(X_t^x),\ \ \ f\in\B_b(\H), t\ge 0, x\in \H,$$  where $\B_b(\H)$ is the set of all bounded measurable   real functions on $\H$.
$P_t$ is called   strong Feller if it sends $\B_b(\H)$ into $C_b(\H),$ the set of all bounded continuous real functions on $\H$.  Our first main result is the following.

\beg{thm}\label{T1.1} Assume  {\bf (a1)}, {\bf (a2)}   and  {\bf (a3)}.   \beg{enumerate}\item[$(1)$] For any $X_0\in\B(\OO\to\H;\F_0)$, the equation $\eqref{E1}$ has a unique mild solution $(X_t)_{t\in [0,\zeta)}$  with life time $\zeta$.
\item[$(2)$] Let $\|Q_t\|_\infty:= \sup_{x\in\H}\|Q_t(x)\|$ be locally bounded in $t\ge 0$.  If  there exist  two positive  increasing functions $\Phi,h: [0,\infty)\times [0,\infty)\to (0,\infty)$ such that  $\int_1^\infty \ff{\d s}{\Phi_t(s)}=\infty $ and
\beq\label{C1} \<(B_t+b_t)(x+y), x\> \le \Phi_t(|x|^2)+ h_t(|y|),\ \ \ x,y\in \H, t\ge 0,\end{equation}   then the mild solution is non-explosive and  $P_t$ is strong Feller for $t>0$.
\end{enumerate}\end{thm}

Without loss of generality, in Theorem \ref{T1.1} one may take   $B=0$ in Theorem \ref{T1.1}. But the situation is different in the next result (Theorem \ref{T1.2}) where the singular part $b$ is bounded in the space variable, so that  the appearance of $B$ allows the whole drift $B_t+b_t$ unbounded and singular.

Comparing with the above mentioned  results of \cite{DF, DR1, DR2}, Theorem \ref{T1.1} contains the following several new points: (1) It works for multiplicative noise; (2) It works for arbitrary starting points and non-H\"older continuous drift; (3) The assertion on the strong Feller property is new, see also  Remark 4.1   for a discussion on Harnack inequalities. Moreover, condition \eqref{C1} is more general than
  $$\<b(x+y),x\>\le C(|x|^2 +1 +\e^{p|y|}),\ \ \ x,y\in\H$$ for some constants $C,p>0$ which is used in  \cite[Theorem 16]{DR2} to enure  the non-explosion of the solution.   See \cite[Remark 17]{DR2} for an explanation on the reasonability of such a condition in infinite dimensions.

The main difficulty in the proof of Theorem \ref{T1.1} comes from the singular drift $b$. To overcome this difficulty, a regularization argument has been introduced in \cite{DF} and  further developed in \cite{DR1, DR2},   to reformulate the mild solution by using a regular functional instead of $b$. This functional is constructed by solving an equation involving in the resolvent associated to the corresponding regular equation, i.e. the equation \eqref{E1} without $b$. Based on such a regularization formulation,   the uniqueness can be proved as in   \cite{FP} where the transport equation for  H\"older continuous vector fields with  a finite-dimensional multiplicative noise is concerned. See also \cite{AG,GP} and references therein for the study of singular SPDEs using regularization by the space-time white noise.

The key point in the proof of Theorem \ref{T1.1}  is to realize the idea of \cite{DF} for the present situation where $Q$ is  non-constant and $b$ is non-H\"older continuous. This is done    by establishing necessary derivative estimates using minimal continuity conditions on $b$.  

\

Next, we consider gradient estimates and Harnack inequalities for the associated Markov semigroup $P_t.$  To this end, we need the following global versions of assumptions {\bf (a2)} and {\bf (a3)}.  For a real function $f$ defined on $[0,T]\times\H$, let  $$\|f\|_{T,\infty}= \sup_{t\in [0,T], x\in\H} |f|(t,x).$$ The same notation applies to $\H$-valued or operator-valued maps, for instance, $\|Q\|_{T,\infty}=\sup_{[0,T]\times \H} \|Q\|.$

\beg{enumerate}\item[{\bf (a2')}]   $B$ and $Q$ satisfy {\bf (a2)}, and moreover
$$  \|b\|_{T,\infty}+\|\nn B\|_{T,\infty} +\|\nn Q\|_{T,\infty} + \|\nn^2 Q\|_{T,\infty}   +\|Q\|_{T,\infty} + \|(QQ^*)^{-1}\|_{T,\infty}\le \Psi(T)$$ holds for some $\Psi\in C([0,\infty))$ and all $T\in [0,\infty)$.
\item[{\bf (a3')}]  $Q$ satisfies \eqref{QQ}. Moreover, for any $T>0$, there exists $\phi\in \D$ such that
\beq\label{BB} |b_t(x)-b_t(y)| \le \phi(|x-y|),\ \ \  t\in [0,T],\, x,y\in\H.  \end{equation}
\end{enumerate}

According to Theorem \ref{T1.1}, under {\bf (a1)}, {\bf (a2')} and {\bf (a3')} the unique mild solution of \eqref{E1} is non-explosive. Let $P_t$ be the associated semigroup.
Gradient estimates and log-Harnack inequalities presented in the next result are new even in finite-dimensions. Note that when $b$ is H\"older continuous and $\H$ is finite-dimensional, the (log) Harnack inequalities have been established recently in \cite{LLW} using the regularization transform of \cite{FP, FGP}. But Theorem \ref{T1.2} also applies  to non-H\"older continuous drifts on infinite-dimensional $\H$.

\beg{thm}\label{T1.2} Assume {\bf (a1)}, {\bf (a2')} and {\bf (a3')}. \beg{enumerate} \item[$(1)$] For any $T>0$ there exists a constant $C(T)>0$ such that
$$ |\nn P_tf|^2 +\ff{P_t f^2- (P_t f)^2} t \le C(T)   P_t|\nn f|^2,\ \ t\in (0,T], f\in C_b^1(\H).$$
\item[$(2)$] There exists a constant $C>0$   such that
\beq\label{G0} |\nn P_t f|^2\le \ff C{t\land 1} \big\{P_t f^2-(P_tf)^2\big\},\ \ t>0, f\in \B_b(\H),\end{equation}
\beq\label{LH0}  P_t\log f(y)\le \log P_t f(x) +\ff{C|x-y|^2}{t\land 1},\ \ t>0, x,y\in \H,  0<f\in \B_b(\H),\end{equation}
\beq\label{H0} P_t f(y)\le P_t f(x) +|x-y|\ss{\ff C{t\land 1} P_t f^2(y)},\ \ x,y\in\H,t>0, 0\le f\in \B_b(\H).\end{equation}
 \end{enumerate} \end{thm}

\paragraph{Remark 1.1.} (1) According to \cite{RW10, WZ14}, the key point in the proof of Theorem \ref{T1.2} is the gradient estimate
$$|\nn P_t f|^2\le C(T) P_t |\nn f|^2.$$ As   the regularization formula \eqref{E4} still contains a non-Lipschitz term $A u_s$, the standard argument in the literature is invalid. Our proof is new in this singular setting (see the proof of Lemma \ref{L6.1}(2)).

(2) In the situation of Theorem \ref{T1.2}, for any $s\ge 0$, let $P_{s,t} f(x)= \E f(X_{s,t}^x), x\in\H, t\ge s, f\in \B_b(\H)$, where $(X_{s,t}^x)_{t\ge s}$ is the unique mild solution to
$$\d X_{s,t}^x = \big\{A X_{s,t}^x +(B_t+b_t)(X_{s,t}^x)\big\}\d t + Q_t (X_{s,t}^x)\d W_t,\ \ t\ge s, X_{s,s}^x =x.$$ Then the assertions in Theorem \ref{T1.2}   hold  for $P_{s,s+t}$ in place of $P_t$ with $C(T)$ and $C$ depending also on $s$. If conditions in {\bf (a2')} and {\bf (a3')} are uniformly in $T$ (i.e. they hold with $T=\infty$ and $[0,\infty)$ in place of $[0,T]$),  then these constants are independent of $s$ and, by the semigroup property, we may take $C(T)= c_1\e^{c_2 T}$ for some constants $c_1,c_2>0$.

\

The inequality in Theorem \ref{T1.2}(1) as well as \eqref{G0} are well known due to Bakry-Emery under a curvature condition, and are easy to check in the regular case (i.e. $b=0$), see e.g. \cite{Bakry, Wb2} and references therein. The log-Harnack inequality \eqref{LH0}
is introduced in  \cite{W10} as a limit version of the dimension-free Harnack inequality founded in \cite{W97}. This inequality has a number of applications. For  instance, it implies that the laws of $X_t^x$ and $X_t^y$ are equivalent and provides pointwise estimates on the Radon-Nikodym derivative; in the time-homogeneous case it implies that the   invariant probability measure $\mu$ (if exists) is unique and has full support on $\H$, the semigroup has positive density (i.e. heat kernel) with respect to $\mu$ (more generally, to an quasi-invariant measure), and it provides heat kernel estimates and entropy-cost inequalities of the semigroup; see, for instance,   \cite[\S 1.4]{Wbook}   for details. Recently, a link of the log-Harnack inequality to the optimal transportation has been presented in \cite{BGL}. Finally, according to  \cite[Proposition 2.3]{ATW14} and \cite[Proposition 1.3]{W14} in a more general framework, the log-Harnack inequality \eqref{LH0}  implies   the gradient estimate \eqref{G0}, while \eqref{G0}  is equivalent to the Harnack type inequality \eqref{H0}.

\

The remainder of the paper is organized as follows. In Section 2, we present some gradient estimates on the  semigroup for the corresponding O-U type  equation, i.e. \eqref{E1} with $B=b=0$. These gradient estimates enable us to prove  the desired regularization representation of the mild solution to \eqref{E1} with non-H\"older drift $b$.   In Section 3, we prove the pathwise uniqueness using the regularization representation  and, in Section 4, we investigate the strong Feller property and discuss Harnack inequalities   for the semigroup. Results  in Sections 3-4 are derived under some global  conditions. Combining these results with   a truncating argument, we prove Theorem \ref{T1.1} in Section 5.  Finally, we prove Theorem \ref{T1.2} in Section 6 by using the regularization representation and finite-dimensional approximations.

\section{Regularization representation of mild solutions}

Since it is easy to construct a weak mild solution of \eqref{E1},   in the spirit of Yamada-Watanabe \cite{YW}     the key point to prove the existence and uniqueness lies in the pathwise uniqueness. To prove the pathwise uniqueness of the mild solution to \eqref{E1}, we aim to construct a  transform $\theta: [0,T]\times \H\to \H$ such that
 \beg{enumerate}
 \item[(a)] For very $t\in [0,T]$, $\theta_t$ is a $C^2$-diffeomorphism on $\H$;
    \item[(b)] If $(X_t)_{t\in [0,T]}$ solves \eqref{E1}, then $\{\theta_t(X_t)\}_{t\in [0,T]}$ solves a regular equation having pathwise uniqueness. \end{enumerate}
  In this way we prove the pathwise uniqueness of \eqref{E1}. For readers' convenience, we briefly explain the idea of the construction of $\theta$ (see also  \cite{DF}).

Write $\theta_t(x)= x+ u_t(x), (t,x)\in [0,T]\times\H$. In order that $\theta_t$ is a $C^2$-diffeomorphism on $\H$, we will take $u_t\in C^1(\H;\H)$ such that $\nn u_t$ is Lipschitz continuous with $\|\nn u_t\|_\infty<1$. By It\^o's formula we have, formally,
\beq\label{**0}\d \theta_t(X_t)= \big\{(\pp_t \theta_t)(X_t)+ (L_t \theta_t)(X_t) \big\}\d t +   (\nn\theta_t)(X_t) \big\{Q_t(X_t)\d W_t +( B_t(X_t)\d t\big\},\end{equation}
where $X_t$ solves \eqref{E1} and
$$L_t:= \ff 1 2\sum_{i,j} \<Q_tQ_t^*e_i,e_j\>\nn_{e_i}\nn_{e_j} + \nn_{A\cdot}+ \nn_{b_t }.$$
To ensure that coefficients in \eqref{**0} are regular as required by point (b),  we set $$\pp_t \theta_t(x) =Ax  -L_t \theta_t(x);$$  i.e. $\pp_t u_t= -L_t u_t-b_t$. In particular, with  $u_T=0$ we have
\beq\label{ABD}u_s=\int_s^TP_{s,t}^0 \{\nn_{b_t} u_t+b_t\}\d t,\ \ \ s\in [0,T],\end{equation} where $\{P_{s,t}^0\}_{0\le s\le t}$ is the semigroup associated to  the O-U type equation
\beq\label{OU}\d Z_{s,t}^x=  A Z_{s,t}^x \d t + Q_t(Z_{s,t}^x)\d W_t,\ \ \ t\ge s, Z_{s,s}^x=x.\end{equation}
It is well known that under assumptions {\bf (a1)} and {\bf (a2')},  the equation \eqref{OU} has a unique mild solution which is non-explosive (see \cite{DZ}). We have
$$P_{s,t}^0f(x)= \E f(Z_{s,t}^x),\ \  \ t\ge s\ge 0, f\in \B_b(\H), x\in\H.$$
To ensure $\|\nn u_s\|_\infty<1$ as required by point (a), instead of \eqref{ABD} we consider
\beq\label{ABC} u_s=\int_s^T\e^{-\ll(t-s)} P_{s,t}^0 \{\nn_{b_t} u_t+b_t\}\d t,\ \ \ s\in [0,T] \end{equation} for large enough $\ll>0$, which    also ensures the desired regularity of the equation \eqref{**0}, see \eqref{O2} and \eqref{O2'} below for details.

\

 To verify the regularity properties of $u_s$ solving \eqref{ABC} for large $\ll>0$, we   first consider  derivative estimates on $P_{s,t}^0$. In the following result,  \eqref{2.3} is more or less standard, but \eqref{2.4} is new.

\beg{lem}\label{L2.1} Assume {\bf (a1)} and {\bf (a2')} with $b=0$. Let $T>0$  be fixed.
\beg{enumerate} \item[$(1)$]  There exists a constant $C>0$ such that for any $f\in \B_b(\H),$
\beq\label{2.3} |\nn P_{s,t}^0 f(x)-\nn P_{s,t}^0 f(y)|\le C\bigg(\ff{|x-y|}{t-s}\land \ff 1{\ss {t-s}}\bigg)\|f\|_\infty,\ \ 0\le s<t\le T,   x,y\in \H.\end{equation}
\item[$(2)$] There exists two constants $c_1,c_2>0$  such that for any increasing  $\phi: [0,\infty)\to [0,\infty)$ with concave $\phi^2$,
\beq\label{2.4}  \|\nn^2 P_{s,t}^0 f\|_{\infty} :=\sup_{x\in\H} \|\nn^2P_{s,t}^0 f(x)\| \le  \ff{c_1 \phi(c_2(t-s)^{\vv/2}) }{t-s}, \ \ 0\le s<t\le T   \end{equation} holds for all $f\in \B_b(\H)$ satisfying
\beq\label{2.5} |f(x)-f(y)|  \le \phi(|x-y|),\ \ \ x,y\in \H,\end{equation} where $\vv\in (0,1)$ is in {\bf (a1)}.
\end{enumerate}\end{lem}

\beg{proof} (1)   We shall make use of the following Bismut formula
\beq\label{C01}  \nn_\eta P_{s,t}^0 f(x)= \E\bigg[\ff{f(Z_{s,t}^x)}{t-s} \int_s^t  \big\<\{Q^*_r(QQ^*)_r^{-1}(Z_{s,r}^x)\nn_\eta Z_{s,r}^{x},\d W_r\big\>_{\bar\H}\bigg]\end{equation} for $x,\eta\in\H, t>s\ge 0,$ and $ f\in\B_b(\H).$ Here, by {\bf (a1)} and {\bf (a2')},    the derivative process  $(\nn_\eta Z_{s,t}^x)_{t\ge s}$ is the unique  mild solution to the linear equation
\beq\label{C0} \d \nn_\eta Z_{s,t}^x=  A \nn_\eta Z_{s,t}^x \d t +\big(\nn_{\nn_\eta Z_{s,t}^x}Q_t\big)(Z_{s,t}^x)\d W_t,\ \ \nn_\eta Z_{s,s}^x=\eta, t\ge s,\end{equation}  so that
\beq\label{C10}\sup_{x\in\H,  0\le s\le t\le T} \E|\nn_\eta Z_{s,t}^x|^2 \le c  |\eta|^2, \ \ \eta\in\H \end{equation} holds for some constant $c>0$; see, for instance, \cite[Remark 9.5]{DZ}.

The formula \eqref{C01} can be easily proved by using the Malliavin calculus. Here, we give a brief proof of this formula for $f\in C_b^1(\H),$ which implies the same formula for $f\in \B_b(\H)$ by an approximation argument. Take
\beq\label{Z0}h_v=\ff 1 {t-s} \int_s^v\{Q^*_r(QQ^*)_r^{-1}\}(Z_{s,r}^x)\nn_\eta Z_{s,r}^x\d r,\ \ v\in [s,t].\end{equation}
In the same manner of \cite[Remark 9.5]{DZ}, but using the Malliavin derivative $D_h$ to replace the directional derivative $\nn_\eta$, we see that
the Malliavin derivative process $(D_h Z_{s,r}^x)_{r\in [s,t]}$ is the unique mild  solution to the equation
$$ \d D_h Z_{s,r}^x=  A D_h Z_{s,r}^x  \d r +\big(\nn_{D_h Z_{s,r}^x}Q_r\big)(Z_{s,r}^x)\d W_r+ Q_r(Z_{s,r}^x)h'_r\d r,\ \
  D_h Z_{s,s}^x=0, r\in [s,t].$$
Combining this with \eqref{C0} and the definition of $h$,
we see that both $(\ff {r-s} {t-s} \nn_\eta Z_{s,r}^x)_{r\in [s,t]}$ and   $(D_h Z_{s,r}^x)_{r\in [s,t]}$ solve   the equation
 $$\d V_r =  A V_r \d r + (\nn_{V_r}Q_r)(Z_{s,r}^x)\d W_r +\ff 1{t-s}  \nn_\eta Z_{s,r}^x\d r,\ \ r\ge s, V_s=0.$$ By the uniqueness of the mild solution to this equation, we obtain $D_h Z_{s,t}^x= \ff {t-s}{t-s}  \nn_\eta Z_{s,t}^x= \nn_\eta Z_{s,t}^x.$ So, by the chain rule and the integration by parts formula in the Malliavin calculus, we arrive at
\beg{equation*}\beg{split} \nn_\eta P_{s,t}^0 f(x) &= \E(\nn_{\nn_\eta Z_{s,t}^x}f)(Z_{s,t}^x)=\E (\nn_{D_h Z_{s,t}^x}f)(Z_{s,t}^x)\\
&=\E D_h (f(Z_{s,t}^x))= \E\bigg[f(Z_{s,t}^x)\int_s^t\<h_r', \d W_r\>_{\bar\H}\bigg].\end{split}\end{equation*} This implies (\ref{C01}).

Now, according to  \eqref{C01}, \eqref{C10} and {\bf (a2')},  there exists a constant $c>0$  such that
\beq\label{C02}  |\nn P_{s,t}^0 f|^2(x)
 \le \ff{c }{t-s} P_{s,t}^0f^2(x),\ \ 0\le s<t\le T, x\in\H, f\in\B_b(\H). \end{equation}  Next, writing $P_{s, t}^0= P_{s, \ff{t+s}2}^0 P_{\ff{t+s}2,t}^0$ by the Markov property,  and applying    \eqref{C01} to $\ff{t+s}2$ and   $P_{\ff{t+s}2,t}^0f$ instead of $t$ and $f$,  we obtain
\beq\label{D*}\nn_\eta P_{s,t}^0 f(x)= \E \bigg[ \ff{(P_{\ff{t+s}2,t}^0f)(Z_{s,\ff{t+s}2}^x)}{(t-s)/2} \int_s^{\ff{t+s}2} \big\<\{Q^*_r(QQ^*)_r^{-1}\}(Z_{s,r}^x) \nn_\eta Z_{s,r}^{x},\d W_r\big\>_{\bar\H}\bigg].\end{equation}
So, for any $\eta'\in\H,$ we can prove
\beq\label{C03} \beg{split} &\ff 1 2 (\nn_{\eta'}\nn_\eta P_{s,t}^0 f)(x)  \\
= \ & \E\bigg[\ff{\big(\nn_{\nn_{\eta'}Z_{s,\ff{t+s}2}^x}P_{\ff{t+s}2,t}^0 f\big)  (Z_{s,\ff{t+s}2}^x) } {t-s}\int_s^{\ff{t+s}2} \big\<\{Q^*_r(QQ^*)_r^{-1}\}(Z_{s,r}^x) \nn_\eta Z_{s,r}^{x},\d W_r\big\>_{\bar\H}\bigg]\\
& +  \E\bigg[\ff{(P_{\ff{t+s}2,t}^0 f)(Z_{s,\ff{t+s}2}^x)}{t-s}\int_s^{\ff{t+s}2}   \big\<\big(\nn_{\nn_{\eta'}Z_{s,r}^x}\{Q^*_r(QQ^*)_r^{-1}\}\big)(Z_{s,r}^x)\nn_\eta Z_{s,r}^{x},\d W_r\big\>_{\bar\H}\bigg]\\
&  +  \E\bigg[\ff{(P_{\ff{t+s}2,t}^0 f)(Z_{s,\ff{t+s}2}^x)}{t-s}\int_s^{\ff{t+s}2}   \big\<\{Q^*_r(QQ^*)_r^{-1}\}(Z_{s,r}^x) \nn_{\eta'}\nn_\eta Z_{s,r}^{x},\d W_r\big\>_{\bar\H}\bigg].\end{split}\end{equation}To verify this formula, we need to apply the dominated convergence theorem.
In the spirit of  \cite[Remark 9.5]{DZ}, {\bf (a1)}, {\bf (a2')} and \eqref{C0}   imply  that  $\ggm_r:= \nn_{\eta'}\nn_\eta Z_{s,r}^x$  is the unique mild solution to   the  equation
$$\d \ggm_r=  A \ggm_r \d r
 +\Big\{\big(\nn_{\ggm_r}Q_r)(Z_{s,r}^x)+\big(\nn_{\nn_{\eta'}Z_{s,r}^x}\nn_{\nn_\eta Z_{s,r}^x}Q_r\big)(Z_{s,r}^x)\Big\}\d W_r,\ \ r\ge s, \ggm_{s}=0,$$ so that
 \beq\label{C04}\sup_{x\in\H, 0\le s\le r\le T} \E |\nn_{\eta'}\nn_\eta Z_{s,r}^x|^2\le c |\eta|^2|\eta'|^2,\ \  \eta,\eta'\in\H \end{equation} holds for some constant $c>0$.    Combining this with \eqref{C02}, \eqref{C10} and {\bf (a2')},  we derive  \eqref{C03}   from \eqref{D*} by using the dominated convergence theorem. Moreover, \eqref{C03}    implies
\beq\label{*3}  |\nn_{\eta'}\nn_\eta P_{s,t}^0 f|^2(x)  \le \ff{c |\eta|^2|\eta'|^2}{(t-s)^2} P_{s,t}^0 f^2(x),
\ \ x,\eta,\eta'\in\H, 0\le s<t\le T, f\in \B_b(\H) \end{equation}  for some constant   $c>0$.    In particular,
$$|\nn P_{s,t}^0f(x)-\nn P_{s,t}^0f(y)|\le \ff {c |x-y|}{t-s}\|f\|_\infty,\ \ x,y\in \H, 0\le s<t\le T, f\in\B_b(\H)$$ holds for some constant $c>0$.  Combining this with \eqref{C02} we prove \eqref{2.3}.

(3)  Applying \eqref{*3} to $\tt f:= f- f(\e^{(t-s)A}x)$ in place of $f$, we obtain
$$|\nn_{\eta'}\nn_\eta P_{s,t}^0 f|^2(x) = |\nn_{\eta'}\nn_\eta P_{s,t}^0 \tt f|^2(x) \le \ff{c |\eta|^2|\eta'|^2}{(t-s)^2} \E \big|f(Z_{s,t}^x)-f(\e^{(t-s)A} x)\big|^2.$$
Since $Z_{s,t}^x- \e^{(t-s)A}x= \int_s^{t} \e^{(t-r)A}Q_r(Z_{s,r}^x)\d W_r=:\bb_{s,t},$ by  \eqref{2.5} and noting that $\phi^2$ is concave and increasing, we obtain
\beq\label{XX} |\nn_{\eta'}\nn_\eta P_{s,t}^0 f|^2(x)
 \le \ff{c |\eta|^2|\eta'|^2}{(t-s)^2} \E \phi^2(|\bb_{s,t}|) \le  \ff{c |\eta|^2|\eta'|^2}{(t-s)^2}   \phi^2\big((\E|\bb_{s,t}|^2)^{\ff 1 2}\big).
  \end{equation}
Since, due to {\bf (a1)},
\beg{equation*}\beg{split} &\E |\bb_{s,t}|^2=  \int_s^{t}\|\e^{(t-r)A} Q_r(Z_{s,r}^x)\|_{HS}^2\d r \\
&\le \|Q\|^2_{T,\infty} \sum_{n=1}^\infty \ff{1-\e^{-2\ll_n (t-s)}}{\ll_n}
\le \|Q\|^2_{T,\infty} \sum_{n=1}^\infty \ff{(2\ll_n (t-s))^\vv}{\ll_n}\le c  (t-s)^\vv\end{split}\end{equation*} holds for $c :=\|Q\|^2_{T,\infty} 2^\vv\sum_{n=1}^\infty\ff 1 {\ll_n^{1-\vv}}<\infty,$  \eqref{2.4} follows from    \eqref{XX}.
\end{proof}

As  a straightforward consequence of   Lemma \ref{L2.1}, we have the following result on the resolvent
$$(R^\ll_{s,t} f)(x) := \int_s^t \e^{-(r-s)\ll}P_{s,r}^0f_r(x) \d r,\ \ \ x\in\H, \ll\ge 0, t \ge s\ge 0, f\in \B_b([0,\infty)\times\H).$$

\beg{lem}\label{L2.2} Assume {\bf (a1)} and {\bf (a2')} with $b=0$. Let $T>0$ be fixed.
\beg{enumerate}\item[$(1)$] There exists a constant $C>0$, such that for any $f\in \B_b([0,T]\times\H),$
$$|\nn (R_{s,t}^\ll f)(x)-\nn (R_{s,t}^\ll f)(y)|\le  C\|f\|_\infty  |x-y| \log \Big(\e+\ff 1 {|x-y|}\Big),\ \ \ll\ge 0, 0\le s\le t \le T.$$
\item[$(2)$] For any $\phi\in\D$, there exists a decreasing   function $\dd_\phi: [0,\infty)\to (0,\infty)$  such that
$$ \lim_{\ll\to \infty} \dd_\phi(\ll)= 0,\ \   \|\nn^2  R_{s,t}^\ll f \|_\infty\le \dd_\phi(\ll),\ \  \ll\ge 0,  0\le s\le t\le T
  $$   for all   $f\in \B_b([0,T]\times\H)$ satisfying \beq\label{D01} |f_t(x)-f_t(y)|  \le \phi(|x-y|),\ \ \ x,y\in \H, t\in [0,T].\end{equation}
\end{enumerate}\end{lem}

\beg{proof}   By Lemma \ref{L2.1}(1) and the definition of $R_{s,t}^\ll f$, there exist  constants $C_1,C_2>0$   such that for  any $ f\in \B_b([0,T]\times \H),$
\beg{equation*}\beg{split} &|\nn (R_{s,t}^\ll f)(x)- \nn (R_{s,t}^\ll f)(y)|\le C_1 \|f\|_\infty \int_s^t \e^{-(r-s)\ll}\Big(\ff 1 {\ss{r-s}}\land \ff{|x-y|}{r-s}\Big) \d r\\
&\le  C_1 \|f\|_\infty \bigg(\int_0^{|x-y|^2\land \e^{-1}} \ff{\d r}{\ss r} +|x-y| \int_{|x-y|^2\land\e^{-1}}^1 \ff{\d r} r +|x-y|\int_1^{T\lor 1} \e^{-\ll r} \d r\bigg)\\
&\le C_2 \|f\|_\infty |x-y|\log \Big(\e + \ff 1 {|x-y|}\Big),\ \ \ x,y\in \H, 0\le s\le t\le T.\end{split}\end{equation*} Then (1) is proved.


Next, since for $\phi\in\D$ we have $\int_0^T \ff {\phi(c_2 s^{\vv/2}) }{s  }\d s<\infty,$    Lemma \ref{L2.1}(2) implies the second assertion for
$$\dd_\phi(\ll):= c_1 \int_0^T\ff{\e^{-\ll s} \phi(c_2 s^{\vv/2})}s\d s \downarrow 0 \ \text{as}\ \ll\uparrow\infty.$$
 \end{proof}

In the next result, we characterize the solution  $u_s$ to \eqref{ABC} which  will be used to formulate the mild solution to \eqref{E1}  (see Proposition \ref{P2.3} below).
To prove the formulation in infinite-dimensions, we shall adopt an approximation argument based on the  second assertion of the following result.

 \beg{lem}\label{L2.3} Assume {\bf (a1)} and {\bf (a2')}, and let $T>0$ be fixed. Then there exists a constant $\ll(T)>0$  such that the following assertions hold.  \beg{enumerate} \item[$(1)$]  For any $\ll\ge \ll(T)$,   the equation $\eqref{ABC}$   has a unique solution
$u\in C([0,T]; C_b^1(\H;\H))$.
\item[$(2)$] Assume that
\beq\label{UC} \lim_{r\downarrow 0} \sup_{|x-y|\le r, t\in [0,T]}|b_t(x)-b_t(y)|=0.\end{equation}     Let $P_{s,t}^{\{n\}}$ be defined as $P_{s,t}^0$ for $Q\circ\pi_n$ in place of $Q$, and let $b^{\{n\}}=  b\circ\pi_n$. Then     for any $\ll\ge \ll(T)$ and $n\ge 1$,     the equation
    \beq\label{D00} u_s^{\{n\}} =\int_s^T \e^{-\ll(t-s)} P_{s,t}^{\{n\}}\big(\nn_{b_t^{\{n\}}}u_t^{\{n\}}+ b_t^{\{n\}}\big)\d t,\ \ s\in [0,T] \end{equation} has a unique  solution $u^{\{n\}}\in C([0,T]; C_b^1(\H;\H))$ such that
\beq\label{0D} \beg{split} & \sup_{n\ge 1}(\|\nn u^{\{n\}}\|_{T,\infty} +\|  u^{\{n\}}\|_{T,\infty})\le \dd(\ll),\\
      &\lim_{n\to\infty} u^{\{n\}} = u,\ \ \ \lim_{n\to\infty} \int_0^T \|\nn u_s -\nn u_s^{\{n\}}\|\d s =0, \end{split} \end{equation}where $\dd$ is a function   such that $\dd(\ll)\to 0$ as $\ll\to\infty$.  If, moreover,  $b$ satisfies  $\eqref{BB}$  for some  $\phi\in\D,$   then
      \beq\label{00D} \sup_{n\ge 1}   \| \nn^2 u ^{\{n\}}\|_{T,\infty}\le \dd_\phi(\ll),\ \ \ll\ge \ll(T) \end{equation} holds for some positive function $\dd_\phi$   such that $\lim_{\ll\to\infty} \dd_\phi(\ll)=0$.  \end{enumerate}\end{lem}

\beg{proof} We first observe that although Lemma \ref{L2.2} is stated for real functions $f$, it works also for $\H$-valued functionals. For instance, if $f\in \B_b([0,T]\times \H;\H)$ satisfies \eqref{D01}, then for any unit $e\in\H$ the real function $\<f,e\>$ satisfies \eqref{D01} as well, so that Lemma \ref{L2.2}(3) implies
$$\sup_{|\eta|\lor |\eta'|\le 1} \|\<\nn_\eta\nn_{\eta'}R_{s,t}^\ll f, e\>\|_\infty  = \sup_{|\eta|,|\eta'| \le 1} \| \nn_\eta\nn_{\eta'}R_{s,t}^\ll \<f, e\>\|_\infty\le \dd_\phi(\ll)|e|,\ \ 0\le s\le t\le T.$$ That is, Lemma \ref{L2.2}(3) works also for $\H$-valued functions.
Below we prove assertions (1) and (2) respectively.

 (1) Let  $\scr H= C([0,T]; C_b^1(\H;\H))$, which is a Banach space under the norm
$$\|u\|_{\scr H} :=\|u\|_{T,\infty}+\|\nn u\|_{T,\infty}= \sup_{t\in [0,T], x\in \H}  |u_t(x)|+ \sup_{t\in [0,T], x\in \H}\|\nn u_t (x)\|,\ \ u\in \scr H.$$
For any $u\in \scr H$, define
$$(\Gamma u)_s(x)= \int_s^{T} \e^{-\ll (t-s)}P_{s,t}^0\big(\nn_{b_t} u_t  +b_{t}\big)(x)\d t, \ \ s\in [0,T].$$ By the fixed-point theorem, it suffices to show that for large enough $\ll>0$, the map $\Gamma$ is   contractive  on $\scr H$. For any $u,\tt u\in\scr \H$, by the definition of $\Gamma$  we have
\beq\label{G1} \|\Gamma u-\Gamma\tt u\|_{T,\infty} \le \int_0^T \e^{-\ll t} \|b\|_{T,\infty} \|\nn u-\nn \tt u\|_{T,\infty}\d t= \ff {\|b\|_{T,\infty}} \ll \|\nn u-\nn\tt u\|_{T,\infty}.\end{equation}
Next, by \eqref{C02} and the definition of $\Gamma$, we have
$$\|\nn(\Gamma u-\Gamma\tt u)\|_{T,\infty}\le C_1\|b\|_{T,\infty}\|\nn u-\nn\tt u\|_{T,\infty} \int_0^T \ff{\e^{-\ll t}}{\ss t} \d t\le \ff{C_2}{\ss\ll} \|\nn u-\nn \tt u\|_{T,\infty}$$ for some constants $C_1,C_2>0$. Combining this with \eqref{G1} we may find $\ll_0(T)>0$  such that   the operator  $\Gamma$ is a contraction operator on $\scr H$ when $\ll\ge \ll_0(T).$

(2) Obviously, if  $(B,b,Q)$ satisfies {\bf (a2')}, so does $(B\circ\pi_n, b\circ\pi_n, Q\circ\pi_n)$ uniformly in $n\ge 1.$ By (1),  $(u^{\{n\}})_{n\ge 1}$ are well defined for $\ll\ge \ll_0(T)$. Due to    \eqref{C02},   there exists a constant $C>0$   such that
$$|\nn u_s^{\{n\}}|\le C \int_s^{T} \ff{\e^{-\ll (t-s)}}{\ss {t-s}}\big(\|\nn u^{\{n\}}\|_{T,\infty} +1\big)\d t,\ \ n\ge 1, s\in [0,T].$$ Taking $\ll_1(T)\ge \ll_0(T)$   such that   $C\int_0^T \ff{\e^{-\ll_1(T) t}}{\ss t}\d t\le\ff 1 2$,   we obtain
$$\|\nn u^{\{n\}}\|_{T,\infty}\le 2C\int_0^T\ff{\e^{-\ll s}}{\ss s}\d s,\ \ \ll\ge \ll_1(T), n\ge 1.$$    Combining this with the definition of $u^{\{n\}}$ we prove
$$\sup_{n\ge 1} \big(\|\nn u^{\{n\}}\|_{T,\infty}+ \|u^{\{n\}}\|_{T,\infty}\big)\le\dd(\ll),\ \ \ll\ge \ll_1(T) $$ for some function  $\dd$   with $\dd(\ll)\to 0$ as $\ll\to\infty$.  Moreover, we have
\beq\label{*Y} \beg{split} u_s-u_s^{\{n\}}= &\int_s^{T} \e^{-\ll (t-s)} P_{s,t}^0 \big\{ \nn_{b_t^{\{n\}}} (u_{t} - u_{t}^{\{n\}}) \big\}\d t\\
&+ \int_s^{T} \e^{-\ll (t-s)} P_{s,t}^0 \big\{\nn_{b_t-b_t^{\{n\}}} u_{t}  +b_{t}-b_{t}^{\{n\}}\big\}\d t\\
&+ \int_s^{T} \e^{-\ll (t-s)} \big(P_{s,t}^0 -P_{s,t}^{\{n\}}\big)\big\{\nn_{b_t^{\{n\}}}u_t^{\{n\}}+b_t^{\{n\}}\big\}\d t.\end{split}\end{equation}
To prove the limits in \eqref{0D}, let
$$g_s(x)=\limsup_{n\to\infty} \|\nn u_s-\nn u_s^{\{n\}}\|(x), \ \ \ s\in [0,T], x\in\H.$$ Then $g\in \B_b([0,T]\times\H)$.
Combining \eqref{*Y} with \eqref{C02} which also holds for $P_{s,t}^{\{n\}}$ uniformly in $n\ge 1$, we obtain
\beq\label{OD1}\beg{split}g_s  \le &C_1 \limsup_{n\to\infty} \int_s^T \ff{\e^{-\ll(t-s)}}{\ss{t-s}}\Big(\big(P_{s,t}^0\|\nn u_t-\nn u_t^{\{n\}}\|^2\big)^{\ff 1 2} + \big(P_{s,t}^0|b_t-b_t^{\{n\}}|^2\big)^{\ff 1 2}\Big)\d t\\
 &+ \lim\sup_{n\to\infty} \int_s^T \e^{-\ll(t-s)}  \big\|\nn\big(P_{s,t}^0 -P_{s,t}^{\{n\}}\big)\big\{\nn_{b_t^{\{n\}}} u_t^{\{n\}}+b_t^{\{n\}}\big\}\big\|\d t \end{split}\end{equation} for some constant $C_1>0$. Obviously, by the dominated convergence theorem we have
 $$\limsup_{n\to\infty} \int_s^T \ff{\e^{-\ll(t-s)}}{\ss{t-s}} \big(P_{s,t}^0\|b_t-b_t^{\{n\}}\|^2\big)^{\ff 1 2} \d t=0.$$
Moreover, by \eqref{UC} and Lemma \ref{L2.2}(1), for every $t\in [0,T]$, the following Lemma \ref{LL} applies to
$$f_n:= \nn_{b_t^{\{n\}}}u_t^{\{n\}}+b_t^{\{n\}},\ \ \ n\ge 1,$$ so that
$$\lim_{n\to \infty} \big\|\nn\big(P_{s,t}^0 -P_{s,t}^{\{n\}}\big)\big\{\nn_{b_t^{\{n\}}} u_t^{\{n\}}+b_t^{\{n\}}\big\}\big\|=0.$$
Combining this with
$$\big\|\nn\big(P_{s,t}^0 -P_{s,t}^{\{n\}}\big)\big\{\nn_{b_t^{\{n\}}}u_t^{\{n\}}+b_t^{\{n\}}\big\}\big\|\le\ff C{\ss{t-s}}$$   for some constant $C>0$
according to \eqref{C02} which also holds for $P_{s,t}^{\{n\}}$ in place of $P_{s,t}^0$, we can apply  the dominated convergence theorem to obtain
$$\limsup_{n\to\infty} \int_s^T \e^{-\ll(t-s)}  \big\|\nn\big(P_{s,t}^0 -P_{s,t}^{\{n\}}\big)\big\{\nn_{b_t^{\{n\}}} u_t^{\{n\}}+b_t^{\{n\}}\big\}\big\|\d t=0.$$ Thus, it follows from \eqref{OD1} that
$$g_s\le C_1 \int_s^T  \ff{\e^{-\ll(t-s)}}{\ss{t-s}} \big(P_{s,t}^0g_t^2\big)^{\ff 1 2}\d t \le \ff{C_2}{\ss\ll} \|g\|_{T,\infty},\ \ s\in [0,T], \ll\ge \ll_1(T)$$ holds for some constant $C_2>0.$ Taking $\ll_2(T)=\ll_1(T) \lor(4 C_2^2)$, we obtain $\|g\|_{T,\infty}\le \ff 1 2\|g\|_{T,\infty}$ for $\ll\ge \ll_2(T).$ Since $g$ is bounded, this implies
$$\lim_{n\to\infty} \|\nn u_s-\nn u_s^{\{n\}}\|=0,\ \ s\in [0,T]$$ provided $\ll\ge \ll_2(T)$.
 Combining this with \eqref{*Y} and Lemma \ref{LL} below, and using again  the dominated convergence theorem, we obtain
 $\limsup_{n\to\infty}|u_s-u_s^{\{n\}}|  =0.$ Therefore, \eqref{0D} holds for $\ll\ge\ll_2(T).$

Finally, let $b$ satisfy  $\eqref{BB}$ for some  $\phi\in\D$. Then, by Lemma \ref{L2.2}(1), there exists a constant $c>0$ such that
the functions $$f_t^{\{n\}}:=  \nn_{b_t^{\{n\}}} u^{\{n\}}_{t}  + b^{\{n\}}_{t},\ \ n\ge 1, t\in [0,T]$$ satisfy $\eqref{2.5}$    with $\tt\phi(s):= c\ss{\phi^2(s)+s}$ in place of $\phi$. Obviously, $\phi\in\D$ implies $\tt\phi\in\D$. So,   \eqref{00D} follows from Lemma \ref{L2.2}(3).  In conclusion, Lemma \ref{L2.3} holds for $\ll(T)=\ll_2(T).$ \end{proof}

\beg{lem}\label{LL}Let $P_{s,t}^{\{n\}}$ be in Lemma $\ref{L2.3}$. For any sequence  $\{f_n\}_{n\ge 1}\subset C_b(\H;\H)$ such that
$$\sup_{n\ge 1} \|f_n\|_\infty<\infty,\ \ \dd_r:= \sup_{|x-y|\le r, n\ge 1} |f_n(x)-f_n(y)|\downarrow 0\ \text{as}\ r\downarrow 0,$$ there holds
$$\lim_{n\to\infty} \Big(\big|P_{s,t}^0f_n-P_{s,t}^{\{n\}} f_n\big|+\big\|\nn P_{s,t}^0 f_n- \nn P_{s,t}^{\{n\}}f_n\big\|\Big)=0,\ \ 0\le s< t.$$\end{lem}
 \beg{proof} Let $(Z_{s,t}^{\{n,x\}})_{t\ge s}$ solve the equation
\beq\label{OU'}\d Z_{s,t}^{\{n,x\}}=  A Z_{s,t}^{\{n,x\}} \d t + Q_t(\pi_n Z_{s,t}^{\{n,x\}})\d W_t,\ \ \ t\ge s, Z_{s,s}^{\{n,x\}}=x.\end{equation} We have $P_{s,t}^{\{n\}}f(x)= \E f(Z_{s,t}^{\{n,x\}}),\ f\in\B_b(\H).$ By {\bf (a2')} it is easy to see that
 \beq\label{WJ2} \lim_{n\to\infty} \E\Big(\big|Z_{s,t}^x- Z_{s,t}^{\{n,x\}}\big|^2+\big\|\nn Z_{s,t}^x- \nn Z_{s,t}^{\{n,x\}}\big\|^2\Big)=0.\end{equation} Then for any $r>0$,
 \beq\label{WJJ}\beg{split} & \limsup_{n\to\infty} \big|P_{s,t}^{\{n\}}f_n(x)- P_{s,t}^0f_n(x)\big|\\
 &\le \dd_r + \limsup_{n\to\infty} \E\big(\big| f_n(Z_{s,t}^x)- f_n(Z_{s,t}^{\{n,x\}})\big|1_{\{|Z_{s,t}^x-Z_{s,t}^{\{n,x\}}|\ge r\}}\big)=\dd_r.\end{split}\end{equation}
 By letting $r\downarrow 0$ we prove $\lim_{n\to\infty}|P_{s,t}^0f_n-P_{s,t}^{\{n\}}f_n|=0.$

 Next, by \eqref{C01} and the corresponding formula for $P_{s,t}^{\{n\}}$, for any $\eta\in\H$ we have
\beg{equation*}\beg{split} &|\nn_\eta P_{s,t}^{\{n\}}f_n(x)- \nn_\eta P_{s,t}^0f_n(x)|\\
&\le  \E\bigg|\ff{f_n(Z_{s,t}^x)- f_n(Z_{s,t}^{\{n,x\}})}{t-s}\int_s^T \big\<Q_r^*(QQ^*)_r^{-1}(Z_{s,r}^x)\nn_\eta Z_{s,r}^x, \d W_r\big\>_{\bar\H}\bigg|\\
&+ \E\bigg|\ff{f_n(Z_{s,t}^{\{n,x\}})}{t-s}\int_s^T \big\<Q_r^*(QQ^*)_r^{-1}(Z_{s,r}^x)\nn_\eta Z_{s,r}^x- Q_r^*(QQ^*)_r^{-1}(\pi_n Z_{s,r}^{\{n,x\}})\nn_\eta Z_{s,r}^{\{n,x\}}, \d W_r\big\>_{\bar\H}\bigg|\\
&=: J_n+J_n'.\end{split}\end{equation*} Similarly to \eqref{WJJ}, we can prove $\lim_{n\to\infty} J_n= 0$ uniformly in $|\eta|\le 1.$   Moreover,   since
 $$\sup_{n\ge 1}\|f_n\|_\infty+ \|Q^*(QQ^*)^{-1}\|_{T,\infty}+ \|\nn Q^*(QQ^*)^{-1}\|_{T,\infty}<\infty,$$  from \eqref{WJ2} we see that $\lim_{n\to\infty} J_n'= 0$ uniformly in $|\eta|\le 1.$    Therefore, $$\lim_{n\to\infty} \|\nn P_{s,t}^0f_n(x) -\nn P_{s,t}^{\{n\}} f_n(x)\|=0.$$ \end{proof}

Finally, we present   the   regularization representation of the mild solution as explained in the beginning of this section.  When $Q$ is constant and $B=0$, this  result is essentially due to \cite{DF, DR1, DR2}. Recall that $Q$ is called cylindrical if there exists $n\ge 1$ such that
$Q(x)=Q(\pi_n x)$ for all $x\in \H$.

\beg{prp}\label{P2.3} Assume {\bf (a1)}, {\bf (a2')}  and either {\bf (a3')} or 
\beg{enumerate}\item[{\bf (a3'')}] $Q$ is cylindrical  and    
$b_s\in C_b(\H;\H)$ for $s\ge 0.$  \end{enumerate}  For any $T>0$, there exists a constant $\ll(T)>0$ such that for any stopping time
  $\tau$,  any adapted continuous process $(X_t)_{t\in [0, \tau\land T]}$    on $\H$  with $\P$-a.s.
\beq\label{E3} X_t= \e^{tA} X_0 +\int_0^t \e^{(t-s)A}\big\{B_s+b_s\big\}(X_s)\d s +\int_0^t\e^{(t-s)A}Q_s(X_s)\d W_s,\ \ t\in [0,\tau\land T],\end{equation}  and   any $\ll\ge \ll(T),$  there holds $ \P$-a.s.
\beq\label{E4}\beg{split}  X_t =  &\  \e^{tA}(X_0+u_0(X_0))-u_t(X_t) +\int_0^t \e^{(t-s)A}\big\{Q_s +(\nn u_s)Q_s\big\}(X_s)\d W_s\\
&\ +\int_0^t \Big\{(\ll-A)\e^{(t-s)A}u_s+\e^{(t-s)A} (B_s+\nn_{B_s}u_s)\Big\}(X_s)\d s ,
 \ \ t\in [0,\tau\land T],\end{split}\end{equation}where  $u$ solves $\eqref{ABC}$, and $(\nn u_s)z:= \nn_z u_s$ for $z\in\H$.
 \end{prp}

\beg{proof}   As in the proof of \cite[Theorem 7]{DR1} (see also the proof of \cite[Theorem 2]{DR2}), we first  make finite-dimensional approximations such that It\^o's formula applies.   For every $n\in\N$, let    $$B_s^{(n)}= \pi_n B_s\circ\pi_n, \  \ b^{(n)}_s  = \pi_n b_s\circ \pi_n,\ \ Q_s^{(n)}= \pi_n Q_s\circ\pi_n,\ \ \ n\ge 1, s\ge 0.$$  Let $(P_{s,t}^{(n)})_{t\ge s\ge 0}$ be the semigroup of the following SDE on $\H_n$ (note that $A=\pi_n A$ holds on $\H_n$):
\beq\label{D*2} \d   Z_{s,t}^{(n,z)} = \big\{A  Z_{s,t}^{(n,z)}\big\} \d t + Q^{(n)}_t (Z_{s,t}^{(n,z)})\d W_t, \ \ Z_{s,s}^{(n,z)}=z\in\H_n, t\ge s. \end{equation} We have
\beq\label{BD}   P_{s,t}^{(n)} f ( z)= \E f(  Z_{s,t}^{(n,z)}),\ \ \ f\in \B_b(\H_n), t\ge s\ge 0, z\in \H_n.\end{equation} It is easy to see that
$(B^{(n)},Q^{(n)})$   satisfies {\bf (a2')} for $\H_n$ in place of $\H$. Let $\ll(T)>0$ be such that assertions in   Lemma \ref{L2.3} hold. Then for any $\ll\ge \ll(T)$,  there exists   unique $u^{(n)} \in C([0,T]; C_b^1(\H_n;\H_n))$ satisfying
\beq\label{E2} u_s^{(n)}=   \int_s^{T} \e^{-\ll (t-s)}P_{s,t}^{(n)}\big( \nn_{b_t^{(n)}} u_{t}^{(n)}   +b_{t}^{(n)}\big)\d t,\ \ \ s\in [0,T]. \end{equation} Let $G^{(n)}_r=  \nn_{b_r^{(n)}} u^{(n)}_r+ b_r^{(n)}, r\ge 0.$ To regularize this functional, we fix $\dd>0$ and let
$$F_{s,r}(z)= P_{s,\dd+r}^{(n)} G_r^{(n)} (\pi_n z), \ \ 0\le s\le r\le T, z\in \H.$$ Then  $F_{s,r}= F_{s,r}\circ\pi_n$   and,  by \eqref{*3} which also holds for $(P_{s,t}^{(n)},\H_n)$ in place of $(P^0_{s,t},\H)$,
\beq\label{N*} \sup_{0\le s\le r\le T} \big\{\|F_{s,r}\|_\infty+ \|\nn F_{s,r}\|_\infty+ \|\nn^2 F_{s,r}\|_\infty\big\}<\infty.\end{equation} So, by \eqref{D*2} and It\^o's formula, for any $0\le  s\le  r\le  T$, we have
\beq\label{IT} \d F_{s,r}(Z_{r,t}^{(n,z)})  =   L^{(n)}_t F_{s,r}(Z_{r,t}^{(n,z)})\d t +\<\nn F_{s,r}(Z_{r,t}^{(n,z)}), Q^{(n)}_t(Z_{r,t}^{(n,z)})\d W_t\>,\ \ t\ge r,\end{equation}  where,  for any second-order differentiable function $F$ on $\H$,
\beg{equation*}\beg{split} L^{(n)}_t F(z)&:= \ff 1 2 \sum_{i,j=1}^n \<Q_tQ^*_te_i,  e_j\>(z) \nn_{e_i}\nn_{e_j} F(z) + \< A z,\nn F(z)\>\\
&= \ff 12\sum_{k\ge 1} (\nn^2_{Q_t(z)\bar e_k} F)(z) + \< A z,\nn F(z)\>,\ \ z\in \H.\end{split}\end{equation*}
 Here, $\nn^2_e:= \nn_e\nn_e$ for $e\in\H$, and $\{\bar e_k\}_{k\ge 1}$ is an orthonormal basis on $\bar \H$. By \eqref{N*},  $Q\in C([0,\infty)\times\H; \scr  L(\bar \H;\H))$ and noting that $F_{s,r}=F_{s,r}\circ\pi_n$, we have $L_\cdot^{(n)} F_{s,r}\in C_b([0,T]\times \H)$ for any $T>0.$  So,
it follows from  \eqref{IT} and the a.e. right-continuity of $Q_tQ_t^*$ that
\beg{equation}\label{NB} \beg{split} &\ff{\d}{\d s} F_{s,r} ( z) := -\lim_{v\downarrow 0} \ff{F_{s-v,r}(z)-F_{s,r}(z)} v=- \lim_{v\downarrow 0} \ff{\E F_{s,r}(Z_{s-v,s}^{(n,z)})-F_{s,r}(z)} v\\
&= -\lim_{v\downarrow 0} \ff 1v\E \int_{s-v}^s(L^{(n)}_t F_{s,r})(Z_{s-v,t}^{(n,z)})\d t
 =- L^{(n)}_s F_{s,r} (z),\ \ r>0, \text{a.e.} \ s\in (0, r]. \end{split}\end{equation} Let
\beq\label{U} u_s^{(n,\dd)}=\int_s^T\e^{-\ll(t-s)}(P_{s,t+\dd}^{(n)}G_t^{(n)})\circ\pi_n\d t= \int_s^T\e^{-\ll(t-s)}F_{s,t}\d t,\ \ s\in [0,T].\end{equation} It follows from \eqref{N*} and \eqref{NB} that
 \beq\label{S1} \beg{split}\pp_s  u_s^{(n,\dd)} &= (\ll- L^{(n)}_s) u_s^{(n,\dd)} -(P_{s,s+\dd}^{(n)}G_s^{(n)})\circ\pi_n\\
 &=(\ll- L^{(n)}_s) u_s^{(n,\dd)} -\big(P_{s,s+\dd}^{(n)}\big\{\nn_{b_s^{(n)}}u_s^{(n)} +b_s^{(n)}\big\}\big)\circ\pi_n.\end{split}\end{equation}    On the other hand, by \eqref{E3}, $X_s^{(n)}:=\pi_n X_s$ solves the following equation on $\H_n$:
$$\d X_s^{(n)} =  A X_s^{(n)}\d s +   \pi_n \big\{B_s+ b_s\big\}(X_s)\d s +  \pi_n Q_s(X_s)\d W_s,\ \ s\in [0,\tau\land T].$$ Then, by $u^{(n,\dd)}= u^{(n,\dd)}\circ\pi_n$ and It\^o's formula,
\beg{equation*}\beg{split} \d  u_s^{(n,\dd)}(X_s^{(n)}) = &  \big(\nn_{Q_s(X_s)\d W_s}  u_s^{(n,\dd)}\big)(X_s^{(n)} )  + (\nn_{b_s(X_s)} u_s^{(n,\dd)})(X_s^{(n)})  \d s\\
&+\big\{(\pp_s u_s^{(n,\dd)})(X_s^{(n)}) + (L^{(n)}_s u_s^{(n,\dd)}+\nn_{\pi_n B_s}u_s^{(n,\dd)})(X_s)\big\} \d s,\ \ \ s\in [0, \tau\land T].\end{split}\end{equation*}
Combining this with \eqref{S1} and noting that $(\nn_{b_s^{(n)}} u_s^{(n)})\circ\pi_n= (\nn_{b_s} u_s^{(n)})\circ\pi_n,$ we obtain
\beg{equation}\label{WW} \beg{split} &(P_{s,s+\dd}^{(n)}b_s^{(n)})(X_s^{(n)})\d s \\
&= \ll  u_s^{(n,\dd)}(X_s^{(n)})\d s - \d  u_s^{(n,\dd)} (X_s^{(n)}) +  \big\{\nn_{b_s(X_s)}  u_s^{(n,\dd)}- P_{s,s+\dd}^{(n)}\nn_{b_s} u_s^{(n)}\big\}(X_s^{(n)} )\d s \\
&\quad +\big(\nn_{Q_s(X_s)\d W_s} u_s^{(n,\dd)}\big) (X_s^{(n)}) +\ff 1 2 \sum_{k\ge 1} \Big(\nn^2_{\{Q_s(X_s)-Q_s(X_s^{(n)})\}\bar e_k} u_s^{(n,\dd)}\Big)(X_s^{(n)} ) \d s\\
&\quad + (\nn_{\pi_nB_s(X_s)} u_s^{(n,\dd)})(X_s^{(n)})\d s,\ \ \ s\in [0,\tau\land T].
 \end{split}\end{equation}  Finally, we complete  proof by using {\bf (a3')} and {\bf (a3'')} respectively.   

(i) Assume {\bf (a3'')}.  Then $(QQ^*)(X_t^{(n)})- (QQ^*)(X_t)=0$ for large $n$, so that \eqref{WW} reduces to $\P$-a.s.
\beg{equation}\label{WW0} \beg{split} &\int_{t_1}^{t_2} (P_{s,s+\dd}^{(n)}b_s^{(n)})(X_s^{(n)})\d s \\
&=  \int_{t_1}^{t_2}\big\{\ll u_s^{(n,\dd)}+\nn_{\pi_n B_s(X_s)} u_s^{(n,\dd)}\big\}(X_s^{(n)})\d s + u_{t_1}^{(n,\dd)}(X_{t_1}^{(n)})-u_{t_2}^{(n,\dd)}(X_{t_2}^{(n)})\\
  &\quad +\int_{t_1}^{t_2} \big\{\nn_{b_s(X_s)}  u_s^{(n,\dd)}- P_{s,s+\dd}^{(n)}\nn_{b_s} u_s^{(n)}\big\}(X_s^{(n)} )\d s +  \int_{t_1}^{t_2}\Big(\nn_{Q_s(X_s)\d W_s} u_s^{(n,\dd)} \Big)(X_s^{(n)})
\end{split}\end{equation} for $0\le t_1\le t_2\le T\land\tau.$ We claim that when $\dd\downarrow 0$ this yields $\P$-a.s.
\beg{equation}\label{WW1} \beg{split}  &\int_{t_1}^{t_2}  b_s^{(n)}(X_s^{(n)})\d s
 =   \int_{t_1}^{t_2}\big\{\ll u_s^{(n)}+\nn_{\pi_n B_s(X_s)+b_s(X_s)-b_s(X_s^{(n)})} u_s^{(n)}\big\}(X_s^{(n)})\d s\\
  &+ u_{t_1}^{(n)}(X_{t_1}^{(n)})-u_{t_2}^{(n)}(X_{t_2}^{(n)})
   + \int_{t_1}^{t_2}\Big(\nn_{Q_s(X_s)\d W_s} u_s^{(n)}\Big) (X_s^{(n)}),\ \ 0\le t_1\le t_2\le \tau\land T.
\end{split}\end{equation}
Indeed, since $\lim_{\dd\downarrow 0} P_{t,t+\dd}^{(n)}f=f$ holds for $t\ge 0$ and $f\in C_b(\H;\H),$ by the boundedness and continuity of $b_s^{(n)}$ and $G_s^{(n)}$, and noting that \eqref{E2} and \eqref{U} imply
\beq\label{WFG} u_s^{(n,\dd)} - u_s^{(n)} = \int_s^T \e^{-\ll(t-s)} \big(P_{s,t}^{(n)} \{P_{t,t+\dd}^{(n)}G_t^{(n)}-G_t^{(n)}\}\big)\circ\pi_n \d t,\end{equation}
we have $\lim_{\dd\downarrow 0} u^{(n,\dd)}= u^{(n)}$ and $\P$-a.s.
$$\lim_{\dd\downarrow 0}\int_{t_1}^{t_2} (P_{s,s+\dd}^{(n)}b_s^{(n)})(X_s^{(n)})\d s= \int_{t_1}^{t_2} b_s^{(n)} (X_s^{(n)})\d s,\ \ 0\le t_1 \le t_2\le T\land\tau.$$
Moreover, combining \eqref{WFG} with \eqref{C02} which also holds for $P_{s,t}^{(n)}$ in place of $P_{s,t}^0$, we obtain
\beq\label{DM} \lim_{\dd\downarrow 0} \|\nn(u_s^{(n,\dd)} -u_s^{(n)})\|\le \lim_{\dd\downarrow 0} \int_s^T \ff c{\ss{t-s}} \ss{P_{s,t}^{(n)} |P_{t,t+\dd}^{(n)} G_t^{(n)} - G_t^{(n)}|^2}\,\d t=0\end{equation} due to the dominated convergence theorem. Thus, $\P$-a.s.
\beg{equation*}\beg{split} &\lim_{\dd\downarrow 0} \int_{t_1}^{t_2}\big\{\nn_{b_s(X_s)} u_s^{(n,\dd)}\big\}(X_s^{(n)})\d s= \int_{t_1}^{t_2}\big\{ \nn_{b_s(X_s)} u_s^{(n)}\big\}(X_s^{(n)})\d s,\\
&\lim_{\dd\downarrow 0} \int_{t_1}^{t_2}\Big(\nn_{Q_s(X_s)\d W_s} u_s^{(n,\dd)} \Big)(X_s^{(n)})=\int_{t_1}^{t_2}\Big(\nn_{Q_s(X_s)\d W_s} u_s^{(n)} \Big)(X_s^{(n)}),\ \ 0\le t_1 \le t_2\le T\land\tau.\end{split}\end{equation*} So, to deduce \eqref{WW1} from \eqref{WW0} with $\dd\downarrow 0$, it remains to prove
$$\lim_{\dd \downarrow 0} \int_{t_1}^{t_2} \big\{\nn_{b_s}  u_s^{(n,\dd)}- P_{s,s+\dd}^{(n)}\nn_{b_s} u_s^{(n)}\big\} \d s=0,\ \ 0\le t_1\le t_2\le T.$$ This follows since by the boundedness of $b$,  the uniform boundedness and continuous of $\nn_{b_s} u_s^{(n)}, G_s^{(n)},$ and  \eqref{DM}, we have
\beg{equation*}\beg{split}& \limsup_{\dd \downarrow 0} \int_{0}^{T} \big|\nn_{b_s}  u_s^{(n,\dd)}- P_{s,s+\dd}^{(n)}\nn_{b_s} u_s^{(n)}\big| \d s\\
&\le \limsup_{\dd \downarrow 0} \int_{0}^{T} \Big(\big|\nn_{b_s}( u_s^{(n,\dd)} -u_s^{(n)})\big|
+ \big| P_{s,s+\dd}^{(n)}\nn_{b_s} u_s^{(n)}-\nn_{b_s} u_s^{(n)}\big| \Big)\d s =0.\end{split}\end{equation*}
Now, writing \eqref{WW1} as
\beg{equation*}\beg{split}  b_t^{(n)} (X_t^{(n)})\d t =& \big\{\ll u_t^{(n)}+ \nn_{\pi_n B_t(X_t)}   u_t^{(n)}\big\} (X_t^{(n)})\d t \\
&- \d  u_t^{(n)}  (X_t^{(n)}) +  (\nn  u_t^{(n)})(X_t^{(n)} ) Q^{(n)}_t(X_t)\d W_t,\ \ 0\le t\le \tau\land T,\end{split}\end{equation*}
we conclude that  $\P$-a.s. for all $t\in [0,T\land\tau]$ and large enough $n$,
\beg{equation*}\beg{split} &X_t-\e^{tA} X_0-\int_0^t \e^{(t-s)A}B_s(X_s)\d s=\int_0^t\e^{(t-s)A}b_s(X_s)\d s +\int_0^t \e^{(t-s)A} Q_s(X_s)\d W_s \\
&=\int_0^t \e^{(t-s)A}  b_s^{(n)}(X_s^{(n)})\d s +\int_0^t \e^{(t-s)A} Q_s(X_s)\d W_s
  +\int_0^t \e^{(t-s)A}\big\{b_s (X_s)- b_s^{(n)}(X_s^{(n)})\big\}\d s\\
 &=   \e^{tA}  u_0^{(n)}(X_0^{(n)})  - u_t^{(n)}(X_t^{(n)})+\int_0^t(\ll-A)  \e^{(t-s)A} u_s^{(n)}(X_s^{(n)})\d s\\
&\quad +\int_0^t \e^{(t-s)A}\Big\{b_s(X_s) - b_s^{(n)}(X_s^{(n)})+\big(\nn_{\pi_n B_s(X_s) }u_s^{(n)}\big)(X_s^{(n)})\Big\} \d s\\
&\quad + \int_0^t \e^{(t-s)A} \Big\{Q_s(X_s)+ (\nn u_s^{(n)})(X_s^{(n)} ) Q^{(n)}_s(X_s)\Big\}\d W_s. \end{split}\end{equation*}
Since $b_s$ and $Q_s$ are bounded and continuous, and    $\|u^{(n)}\|_\infty+\|\nn u^{(n)}\|_\infty$  is bounded in $n$ by Lemma \ref{L2.3}(2),  with $n\to \infty$ this implies \eqref{E4} provided
\beq\label{D*0}\lim_{n\to\infty} u^{(n)}\circ \pi_n  = u,\ \  \ \lim_{n\to\infty} \int_0^T\|\nn u_s-\nn u^{(n)}_s\circ\pi_n\| \d s =0,\end{equation}where the first limit implies $\int_0^t (\ll-A)\e^{A(t-s)} u_s^{(n)} (X_s^{(n)})\d s\to \int_0^t (\ll-A)\e^{A(t-s)} u_s (X_s)\d s$  weakly in $\H$ as $n\to\infty$.
To prove \eqref{D*0} using Lemma \ref{L2.3}(2), let $(Z_{s,t}^{\{n,z\}})_{t\ge s}$ solve the equation \eqref{OU'} for $z$ in place of $x$. Since $\pi_n A=A$ holds on $\H_n$, we see that $\pi_n Z_{s,t}^{\{n,z\}}$ solves \eqref{D*2} for $\pi_n z$ in place of $z$. Thus, $\pi_n Z_{s,t}^{\{n,z\}}= Z_{s,t}^{\{n,\pi_n z\}}$, so that
$$P_{s,t}^{(n)} f(\pi_n z) = P_{s,t}^{\{n\}}(f\circ \pi_n)(z),\ \ z\in\H, f\in \B_b(\H).$$ Combining this with \eqref{E2} and $b^{(n)}\circ \pi_n= b^{(n)}$, we conclude that $u^{\{n\}}:= u^{(n)}\circ\pi_n$ solves \eqref{D00}.
  Therefore,   \eqref{D*0}  follows from   Lemma \ref{L2.3}(2).

(ii) Assume {\bf (a3')}. Then \eqref{QQ} and  \eqref{BB} hold  for some $\phi\in\D$. By Lemma \ref{L2.3}(2),   $\|\nn u^{(n)}\|_\infty+\|\nn^2 u^{(n)}\|_\infty$ is bounded in $n\ge 1$. Since
$$u_s^{(n,\dd)}= \int_s^{T}\e^{-\ll(t-s)}( P_{s,t+\dd}^{(n)} G_t^{(n)})\circ\pi_n\d t,$$ and as explained in the proof of \eqref{00D} that $f_t:=1_{[s,T]}(t)(P_{s,t+\dd}^{(n)}G^{(n)}_t)\circ\pi_n$ satisfies \eqref{2.5} for some $\tt\phi\in\D$, we have $\sup_{n\ge 1, \dd\in (0,1)} \|\nn^2 u^{(n,\dd)}\|_\infty<\infty$ according to Lemma \ref{L2.3}(2).
 Combining this with \eqref{QQ}, we obtain
$$\limsup_{n\to\infty}\limsup_{\dd\downarrow 0} \int_0^{\tau\land T } \sum_{k\ge 1} \big|\nn^2_{\{Q_s(X_s^{(n)})-Q_s(X_s)\}\bar e_k}u_s^{(n,\dd)}(X_s^{(n)})\big|  \d s=0.$$ Therefore, repeating the argument in case (i) we prove \eqref{E4}.
\end{proof}

\section{Pathwise uniqueness}

In this section, we   prove the pathwise uniqueness of mild solutions  under {\bf (a1)}, {\bf (a2')}, and either {\bf (a3')} or the following stronger version  of {\bf (a3)}. \beg{enumerate}\item[{\bf (a3')}]  $Q$ is cylindrical, i.e. $Q=Q\circ\pi_n$; $b\in \B_b( [0,\infty)\times\H;\H)$  such that \eqref{BB} holds for some $\phi\in \D_0.$
 \end{enumerate}

 \beg{prp}\label{P3.1} Assume {\bf (a1)}, {\bf (a2')}   and  {\bf (a3')}.
 Let   $(X_t)_{t\ge 0}, (Y_t)_{t\ge 0}$ be two adapted continuous process on $\H$ with $X_0=Y_0.$ For any $n\ge 1$, let
 $$\tau_n^X=n\land \inf\{t\ge 0: |X_t|\ge n\},\ \ \  \tau_n^Y=n\land \inf\{t\ge 0: |Y_t|\ge n\}.$$ If $\P$-a.s. for all $t\in [0,\tau_n^X\land \tau_n^Y]$ there holds
\beg{equation*} \beg{split} &X_t= \e^{tA}x + \int_0^t \e^{(t-s)A}(B_s+ b_s)(X_s) \d s + \int_0^t \e^{(t-s)A} Q_s(X_s)\d W_s,\\
 & Y_t= \e^{tA}x + \int_0^t \e^{(t-s)A} (B_s+b_s)(Y_s) \d s + \int_0^t \e^{(t-s)A} Q_s(Y_s)\d W_s,\end{split}\end{equation*} then $\P$-a.s. $X_t=Y_t$ for all $t\in [0, \tau_n^X\land\tau_n^Y].$ In particular, $\P$-a.s. $\tau_n^X=\tau_n^Y.$ \end{prp}

 \beg{proof} For any $m\ge 1$, let $$\tau_m=\tau_n^X\land\tau_n^Y\land\inf\{t\ge 0: |X_t-Y_t|\ge m\}.$$  It suffices to prove that for any $T>0$ and $m\ge 1$,
 \beq\label{A1} \int_0^T\E\big[1_{\{s<\tau_m\}}|X_{s}-Y_{s}|^2\big] \d s =0.\end{equation}  Let $\ll>0$ be large enough such that assertions in Lemma \ref{L2.3}  and Proposition \ref{P2.3} hold. By \eqref{E4} for $\tau=\tau_m$, we have $\P$-a.s.
 \beq\label{A2} \beg{split} &X_t-Y_t=   u_t(Y_t)-u_t(X_t) + \int_0^t (\ll-A)\e^{(t-s)A} \big(u_s(X_s)-u_s(Y_s)\big) \d s \\
 & +\int_0^t\e^{(t-s)A}  \big\{(B_s+\nn_{B_s} u_s)(X_s)-(B_s+\nn_{B_s} u_s)(Y_s)\big\}\d s\\
 &+ \int_0^t\e^{(t-s)A} \big(\nn u_s(X_s)-\nn u_s(Y_s)\big) Q_s(X_s)\d W_s\\
 &+ \int_0^t\e^{(t-s)A} \big(\nn u_s(Y_s)+I\big) \big(Q_s(X_s)-Q_s(Y_s)\big)\d W_s,\ \ t\in [0,\tau_m\land T].\end{split}\end{equation} Since $b$ and $u$ are bounded on $[0,T]\times \H$, by \eqref{ABC} and \eqref{C02} we may find a constant $C>0$ such that
 \beq\label{A} \|\nn u_t\|_\infty \le C\int_0^T\ff{\e^{-\ll s}}{\ss s} \d s\le \ff 1 5,\ \ t\in [0,T]\end{equation} for large $\ll>0$. Combining this with \eqref{A2} we obtain $\P$-a.s.
 \beg{equation}\label{WZ1}\beg{split} & |X_t-Y_t|\le  \ff{5} 4 \bigg|\int_0^t (\ll-A)\e^{(t-s)A} \big(u_s(X_s)-u_s(Y_s)\big) \d s\bigg|\\
   &+\ff 54 \bigg| \int_0^t\e^{(t-s)A}  \big(B_s(X_s)-B_s(Y_s)+\nn_{B_s(X_s)-B_s(Y_s)} u_s(X_s)\big)\d s\bigg|\\
 &+ \ff 54 \bigg| \int_0^t\e^{(t-s)A}\big(\nn u_s(Y_s)+I\big) \big(Q_s(X_s)-Q_s(Y_s)\big)\d W_s \bigg|\\
 &+ \ff 54 \bigg| \int_0^t\e^{(t-s)A}  \big(\nn_{B_s(Y_s)}u_s(X_s) - \nn_{B_s(Y_s)}u_s(Y_s)\big)\d s\bigg|\\
  & +\ff 54 \bigg| \int_0^t\e^{(t-s)A} \big(\nn u_s(X_s)-\nn u_s(Y_s)\big) Q_s(X_s)\d W_s \bigg|,\ \  t\in [0,\tau_m\land T].\end{split}\end{equation}  Moreover, by {\bf (a1)} there exists some function $\vv(\ll)\downarrow 0$ as $\ll\uparrow \infty$ such that
\beg{equation}\label{WZ2}\beg{split} &\int_0^r \e^{-2\ll t} \E \bigg|1_{\{t<\tau_m\}}\int_0^{t} \e^{(t-s)A} \big(\nn u_s(X_s)-\nn u_s(Y_s)\big) Q_s(X_s)  \d W_s\bigg|^2\d t \\
  &\le \|Q\|_{T,\infty}^2 \int_0^r \e^{-2\ll t} \d t \int_0^t  \|\e^{A(t-s)}\|_{HS}^2  \E \big[1_{\{s<\tau_m\}}\big|\nn u_{s}(X_s)-\nn u_s(Y_s)\big|^2\big]\d s\\
  &= \|Q\|_{T,\infty}^2\int_0^r \e^{-2\ll s}\E \big[1_{\{s<\tau_m\}}\big|\nn u_s(X_s)-\nn u_{s}(Y_{s})\big|^2\big]\d s\int_s^r \|\e^{A(t-s)}\|_{HS}^2 \e^{-2\ll(t-s)}\d t\\
  &\le \vv(\ll) \int_0^r \e^{-2\ll s}\E \big[1_{\{s<\tau_m\}}\big|\nn u_{s}(X_{s})-\nn u_{s}(Y_{s})\big|^2\big]\d s,\ \ r\in [0,T].\end{split}\end{equation} Similarly, since \eqref{A} and {\bf (a2')} imply $\|\nn u\|_{T,\infty}+ \|\nn Q\|_{T,\infty}<\infty$,  
  \beg{equation*} \beg{split} &\int_0^r \e^{-2\ll t} \E \bigg|1_{\{t<\tau_m\}}\int_0^{t} \e^{(t-s)A} \big(\nn u_s(Y_s)+I\big)\big(Q_s(X_s)-Q_s(Y_s)\big)   \d W_s\bigg|^2\d t \\
    &\le \vv(\ll) \int_0^r \e^{-2\ll s}\E \big[1_{\{s<\tau_m\}}\big| X_{s} - Y_{s} \big|^2\big]\d s,\ \ r\in [0,T] \end{split}\end{equation*}holds for the same type $\vv(\ll)$. 
  Combining this with \eqref{WZ1} and \eqref{WZ2}, and  using {\bf (a2')} and \eqref{A}, we may find a   constant $C_0>0$    such that for large enough $\ll>0$ 
 \beq\label{A3} \beg{split} &\eta_r: = \int_0^r \e^{-2\ll t} \E \big[1_{\{t<\tau_m\}}|X_{t} - Y_{t}|^2\big]\d t \\
 &\le \ff{25}2\int_0^r \e^{- 2\ll t} \E  \bigg|1_{\{t<\tau_m\}}\int_0^{t} (\ll-A)\e^{(t-s)A} \big(u_s(X_s)-u_s(Y_s)\big) \d s\bigg|^2\d t+C_0\int_0^r \eta_t\d t\\
 &\quad + C_0 \int_0^r  \e^{-2\ll s}  \E\big[1_{\{s<\tau_m\}}\|\nn u_{s}(X_{ s})-\nn u_{ s}(Y_{ s})\|^2\big]\d s , \ \ r\in [0,T].\end{split} \end{equation}
    Since
 \beg{equation*}\beg{split} I_t &:=  \e^{-2\ll t}  \bigg|1_{\{t<\tau_m\}}\int_0^{t} (\ll-A)\e^{(t-s)A} \big(u_s(X_s)-u_s(Y_s)\big) \d s\bigg|^2\\
 &= \sum_{i=1}^\infty \e^{-2\ll t} \bigg|(\ll+\ll_i)\int_0^{t\land\tau_m} \e^{-(t\land\tau_m-s)\ll_i}\<u_s(X_s)-u_s(Y_s), e_i\>\d s\bigg|^2\\
 &\le \sum_{i=1}^\infty \bigg(\int_0^{t} (\ll+\ll_i) \e^{-(t-s)(\ll+\ll_i)}\d s\bigg)\\
 &\qquad \times \int_0^{t\land\tau_m} (\ll+\ll_i) \e^{-(t-s)(\ll+\ll_i)-2\ll s}\<u_s(X_s)-u_s(Y_s), e_i\>^2 \d s\\
 &\le \sum_{i=1}^\infty \int_0^{t\land\tau_m} (\ll+\ll_i) \e^{-(t\land\tau_m-s)(\ll+\ll_i)-2\ll s}\<u_s(X_s)-u_s(Y_s), e_i\>^2 \d s,\end{split}\end{equation*} it follows from  \eqref{A} that
 \beq\label{A4} \beg{split}& \E\int_0^r I_t\d t \\
 &\le \sum_{i=1}^\infty (\ll+\ll_i) \E\int_0^r\d t  \int_0^{t\land\tau_m}   \e^{-(t-s)(\ll+\ll_i)-2\ll s}\<u_s(X_s)-u_s(Y_s), e_i\>^2 \d s\\
 &= \sum_{i=1}^\infty \E\int_0^{r\land\tau_m} \e^{-2\ll s} \<u_s(X_s)-u_s(Y_s), e_i\>^2 \d s \int_s^{r} (\ll+\ll_i) \e^{-(t-s)(\ll+\ll_i)}\d t\\
 &\le  \E\int_0^r \e^{-2\ll s} 1_{\{s<\tau_m\}}|u_{s}(X_{s})- u_{s}(Y_{s})|^2 \d s \\
 &\le \ff 1 {25} \int_0^r\e^{-2\ll s}\E\big[1_{\{s<\tau_m\}}|X_{s}-Y_{s}|^2\big]\d s=\ff 1 {25} \eta_r.\end{split}\end{equation}

 Next, by the boundedness of $b$  and Lemma \ref{L2.3}(1), we have $\|\nn_b u+b\|_{T,\infty}<\infty.$   So, according to  Lemma \ref{L2.2}(1) and \eqref{ABC}, there exists a constant $C_1>0$ such that
 $$\|\nn u_s(x)-\nn u_s(y)\|\le C_1 |x-y|\log(\e +|x-y|^{-1}),\ \ s\in [0,T], x,y\in \H.$$
 If \eqref{BB} holds for some $\phi\in \D$, then by Lemma \ref{L2.2}(1) and  $\|b\|_{T,\infty}+\|\nn u\|_{T,\infty}<\infty,$ we conclude that
 $$f_t:=  \nn_{b_t}  u_{t} + b_{t},\ \ \ t\in [0,T]$$ satisfies \eqref{D01} with $\phi$ replaced by $\tt\phi(s):= c\ss{\phi(s)^2+s}$, which is in $\D$ as well.   Therefore,  by {\bf (a1)} and Lemma \ref{L2.3}(2),  when $\ll$ is large enough we have
$$C_0 \int_0^r     \E\big[1_{\{s<\tau_m\}}\|\nn u_{ s}(X_{ s})-\nn u_{ s}(Y_{ s})\|^2\big]\d s  \le \ff 1 4 \eta_r,\ \ \ r\in [0,T].$$
Substituting this and   \eqref{A4} into \eqref{A3}, we arrive at
 $$\eta_r\le \ff 3 4 \eta_r  +C_0\int_0^r \eta_t\d t,\ \ \ r\in [0,T].$$ Hence,
\beq\label{FJ}\eta_r\le 4 C_0 \int_0^r \eta_t  \d t,\ \ \ r\in [0,T].\end{equation}  By the Gronwall inequality we obtain $\eta_T=0,$  which is equivalent to  the desired \eqref{A1}.
 \end{proof}

\section{Strong Feller property and Harnack inequality}

In this section, we investigate the strong Feller property  and discuss Harnack inequalities of the   semigroup associated to the equation \eqref{E1}.

 \beg{prp}\label{P3.2} Let $B=0$, $b_t\in C_b(\H;\H),$ and $Q_t\in C_b^1(\H;\scr L(\bar\H;\H))$    for every $t\ge 0$. Assume
\beq\label{HH}  \|b\|_{T,\infty}+ \|Q\|_{T,\infty}+\|\nn Q\|_{T,\infty}+ \|(QQ^*)^{-1}\|_{T,\infty}<\infty,\ \  T>0.\end{equation}      If,  for any $x\in \H$ and any cylindrical Brownian motion $(W_t)_{t\ge 0},$  the equation $\eqref{E1}$ has a unique mild solution, then the associated Markov semigroup $P_t$ is strong Feller for $t>0$.\end{prp}

 \beg{proof}   For fixed $z\in\H, T>0$ and $f\in\B_b(\H)$, we intend to prove
 \beq\label{ST} \lim_{x\to z} P_Tf(x)= P_T f(z).\end{equation}
To this end, we   formulate $P_T$ using the mild solution to the regular equation
$$\d Z_t^x = A Z_t^x \d t +Q_t(Z_t^x)\d W_t,\ \ Z_0^x=x.$$    More precisely, we have
\beg{equation*}\beg{split}Z_t^x&:= \e^{tA}x  +\int_0^t \e^{(t-s)A} Q_s(Z_s^x)\d W_s\\
&= \e^{tA} x  +\int_0^t \e^{(t-s)A}  b_s  (Z_s^x)\d s +\int_0^t\e^{(t-s)A} Q_s(Z_s^x)\d   W_s^x, \ \ t\in [0,T],\end{split}\end{equation*}
where $$ W_t^x:= W_t-\int_0^t \{Q^*_s(QQ^*)^{-1}_s b_s\}(Z_s^x))\d s,\ \ t\in [0,T].$$  By the Girsanove theorem,
$(W_t^x)_{t\in [0,T]}$ is a cylindrical Brownian motion on $\bar\H$ under probability $\d\Q^x:= R^x_T\d\P$, where
$$R^x_T:= \exp\bigg[\int_0^T \Big\<\{Q^*_s(QQ^*)^{-1}_sb_s\}(Z_s^x),\d W_s\Big\>_{\bar\H}-\ff 1 2 \int_0^T \Big|\{Q^*_s(QQ^*)_s^{-1}b_s\}(Z_s^x)\Big|^2_{\bar\H}\d s\bigg].$$
Then $(Z_t^x, W_t^x)_{t\in [0,T]}$ is a weak mild solution to \eqref{E1},  so that
\beq\label{ST1} P_Tf(x)= \E\big[f(Z_T^x)R^x_T\big],\ \ \ x\in\H.\end{equation}  By the boundedness and continuity of $Q^*_s(QQ^*)^{-1}_sb_s,$ and noting that $Z_t^x$ is continuous in $x$, we conclude that
\beq\label{Z1} \lim_{x\to z} \big|P_Tf(x)- \E\big[f(Z_T^x)R^z\big]\big| \le \|f\|_\infty \lim_{x\to z} \E|R^x_T-R_T^z|=0.\end{equation}
Next, to prove the continuity of $\E\big[f(Z_T^x)R^z\big]$ in $x$, we approximate $b$ by $C_b^1$ maps such that Malliavin calculus can be applied.
Since $b$ is bounded and continuous in the space variable, we may find a sequence $\{b^{(n)}\}_{n\ge 1}\subset \B_b([0,T]\times\H)$ such that $b_s^{(n)}\in C_b^1(\H;\H)$ for $s\in [0,T]$ with $\|\nn b^{(n)}\|_{T,\infty}<\infty$ for every $n\ge 1$,  $\sup_{n\ge 1}\|b^{(n)}\|_{T,\infty}<\infty,$   and $\lim_{n\to \infty} b^{(n)}=b$ holds on $[0,T]\times\H$. Let
$$R^x_{T,n}:= \exp\bigg[\int_0^T \<\{Q^*_s(QQ^*)_s^{-1}b_s^{(n)}\}(Z_s^x),\d W_s\>_{\bar\H}-\ff 1 2 \int_0^T |\{Q^*_s(QQ^*)_s^{-1}b_s^{(n)}\}(Z_s^x)|^2_{\bar\H}\d s\bigg],\ \ n\ge 1.$$ It is easy to see that $R_{T,n}^z$ is Malliavin differentiable   and
\beq\label{Z2} \lim_{n\to\infty} \sup_{x\in \H} \Big|\E\big[f(Z_T^x)R^z_T\big]- \E\big[f(Z_T^x)R_{T,n}^z\big]\Big| \le \|f\|_\infty\lim_{n\to\infty} \E |R^z_T-R_{T,n}^z|=0.\end{equation}
Now, for $\eta\in\H$, let $h$ be in \eqref{Z0} with $s=0$ and $t=T$  such that $\nn_\eta Z_T^x= D_h Z_T^x$ according to the proof of Lemma \ref{L2.1}(1). Then,    for $\nn_\eta$ being taken with respect to the variable $x$, it follows from the integration by parts formula that for any $f\in C_b^1(\H)$,
\beg{equation*}\beg{split} \nn_\eta\E\big[f(Z_T^x)R^z_{T,n}\big] &= \E\big[(\nn_{\nn_\eta Z_T^x} f)(Z_T^x)R^z_{T,n}\big]= \E\big[D_h\{f(Z_T^x)\}R^z_{T,n}\big]\\
&= \E\big[D_h\{f(Z_T^x) R^z_{T,n}\}\big]- \E\big[ f(Z_T^x) D_hR^z_{T,n}\big]\\
&= \E\bigg\{f(Z_T^x)\bigg(R_{T,n}^z\int_0^T\<h_t',\d W_t\>_{\bar\H} -   D_hR^z_{T,n}\bigg)\bigg\},\ \ f\in C_b^1(\H;\H).\end{split}\end{equation*} Up to an approximation argument this implies that $|\nn_\eta \E [f(Z_T^x)R_{T,n}^z]|<\infty$ for any $f\in\B_b(\H)$, where the derivative $\nn_\eta$ is taken with respect to $x$. In particular, $\E [f(Z_T^x)R_{T,n}^z]$ is continuous in $x$. Combining this with \eqref{Z1} and \eqref{Z2}, we prove \eqref{ST}. \end{proof}

\paragraph{Remark 4.1.} Using  coupling by change of measures   as in \cite[Chapter 4]{Wbook}, in the situation of Proposition \ref{P3.2} we may derive  the dimension-free  Harnack inequality   in the sense of \cite{W97}.   Here, instead of repeating the coupling arguments therein,  we intend to show that   \eqref{ST1} together with known Harnack inequalities of $P_T^0$ implies  the corresponding inequalities for $P_T$.    For instance, when $Q_t(x)=Q_t$ does not depend on $x$, by \cite[Theorem 3.2.1]{Wbook} for $K=0$ and $ \ll_T:= \sup_{t\in [0,T]}\|Q_t^*(Q_tQ_t^*)^{-1}\|^2$, the Harnack inequality
\beq\label{H} (P_T^0f(y))^p \le P_T^0 f^p(x) \exp\Big[\ff {p|x-y|^2}{2\ll_T (p-1)}\Big],\ \ p>1, x,y\in \H, f\in \B_+(\H)\end{equation} holds, where $\B_+(\H)$ is the set of all positive measurable functions on $\H$. On the other hand, by \eqref{ST1} and H\"older's inequality, we have
\beg{equation*}\beg{split} &(P_T f(x))^p\le  P_T^0f^p(x)   \big(\E (R^x_T)^{\ff p{1-p}}\big)^{p-1},\\
&(P_T^0 f(x))^p \le \big(\E f^p(Z_T^x)R^x\big)  \big(\E (R^x_T)^{\ff 1{1-p}}\big)^{p-1} =  P_T^pf(x)  \big(\E (R^x)^{\ff 1 {1-p}}\big)^{p-1},\ \ p>1.\end{split} \end{equation*} Combining these with \eqref{H} we obtain
\beg{equation} \label{H2}\beg{split} (P_T f(x))^{p^3} &\le (P_T^0 f^{p}(x))^{p^2}\big(\E (R^x_T)^{\ff p{1-p}}\big)^{p^2(p-1)}\\
&\le (P_T^0 f^{p^2} (y))^p  \big(\E (R^x_T)^{\ff p{1-p}}\big)^{p^2(p-1)} \exp\Big[\ff {p^2|x-y|^2}{2\ll_T (p-1)}\Big]\\
&\le (P_Tf^{p^3})(y)\big(\E (R_T^y)^{\ff 1 {1-p}}\big)^{p-1} \big(\E (R^x_T)^{\ff p{1-p}}\big)^{p^2(p-1)} \exp\Big[\ff {p^2|x-y|^2}{2\ll_T (p-1)}\Big]\end{split}\end{equation}
 for any $p>1$ and $f\in \B^+(\H).$    When the noise is multiplicative, we may derive the Harnack inequality from \cite[Theorem 3.4.1(2)]{Wbook} for large enough $p>1$.

 Comparing with \eqref{H}, the Harnack inequality included in \eqref{H2} is worse for short distance since
   $$\lim_{y\to x} \big(\E (R_T^y)^{\ff 1 {1-p}}\big)^{p-1} \big(\E (R^x_T)^{\ff p{1-p}}\big)^{p^2(p-1)} \exp\Big[\ff {p^2|x-y|^2}{2\ll_T (p-1)}\Big]>1.$$ In particular, it does not imply the strong  Feller property as \eqref{H} does.
   See Section 6 for the study of the log-Harnack inequality which is sharp for short distance as in the regular case.

\section{Proof of Theorem \ref{T1.1} }

Throughout this section, we assume {\bf (a1)}, {\bf (a2)}  and either {\bf (a3)}. Using $b$ to replace $b+B$, we may and do assume that $B=0.$

  (a) We first assume further that  {\bf (a2')}  and  {\bf (a3')} hold.   In this case, for any  $X_0\in\B(\OO\to\H;\F_0)$,  the equation \eqref{E1} has a weak mild solution as shown in the proof of Proposition \ref{P3.2} for $X_0$ in place of $x$. On the other hand, by Proposition \ref{P3.1} we have the pathwise uniqueness of the mild solution. So, by the Yamada-Watanabe principle \cite{YW}
  (see \cite[Theorem 2]{OD04} or \cite{LR13} for the result in infinite dimensions), the equation \eqref{E1} with $B=0$  has a unique mild solution. Moreover, in this case the solution is non-explosive.

(b) In general,  take $\psi\in C_b^\infty ([0,\infty))$ such that $0\le \psi\le 1, \psi(r)=1$ for $r\in [0,1]$ and $\psi(r)=0$ for $r\ge 2.$  For any $m\ge 1, t\ge 0$ and $z\in\H$,  let
$$  b_t^{[m]}(z)= b_{t\land m}(z) \psi(|z|/m),$$ and
$$Q_t^{[m]}(z) = \beg{cases} Q_t\big(\psi(|\pi_n z|/m)z\big)=Q_t\big(\psi(|\pi_n z|/m)\pi_n z\big),&\text{if}\ Q=Q\circ\pi_n\ \text{for\ some\ }n\ge 1,\\
 Q_t(\psi(|z|/m)z),&\text{otherwise.}\end{cases} $$  By {\bf (a2)} we see that {\bf (a2')} holds for $B=0$ and $Q^{[m]}$   in place of $Q$. Moreover,  {\bf (a3)}   implies that $(Q^{[m]},  b^{[m]})$ satisfies {\bf (a3')}.    Then by (a), \eqref{E1} for $B=0$ and $(b^{[m]}, Q^{[m]})$ in place
of $(b,Q)$ has a unique mild solution $X_t^{(m)}$ starting at $X_0$ which is non-explosive. Let
$$\tau_n=n\land \inf\{t\ge 0: |X_t^{(n)}|\ge n\},\ \ \ n\ge 1.$$   Since   $b_s^{[m]}(z)=b_s(z)$ and $ Q_s^{[m]}(z)=Q_s(z)$ hold for $s\le m$ and $|z|\le m,$   by Proposition \ref{P3.1}, for any $n,m\ge 1$ we have
$X_t^{(n)}=X_t^{(m)}$ for $t\in [0, \tau_n\land\tau_m]$. In particular, $\tau_m$ is increasing in $m$.  Let $\zeta=\lim_{m\to\infty} \tau_m$ and
$$X_t= \sum_{m=1}^\infty 1_{[\tau_{m-1}, \tau_m)}(t) X_t^{(m)},\ \ \ \tau_0:=0, t\in [0,\zeta).$$ Then it is easy to see that $(X_t^x)_{t\in [0,\zeta)}$ is a mild solution to \eqref{E1} for $B=0$  with life time $\zeta$  and,  due to Proposition \ref{P3.1},  the mild solution is unique. We prove Theorem \ref{T1.1}(1) for $B=0$.

(c) Let $\|Q\|_{T,\infty}<\infty$ for $T>0$, and let \eqref{C1} hold for some positive increasing $\Phi, h$ such that $\int_1^\infty \ff{\d s}{\Phi_t(s)} =\infty, t\ge 0.$ Let $(X_t)_{t\in [0,\zeta)}$ be a mild solution to \eqref{E1} for $B=0$ with life time $\zeta$. Let
$\xi_t=\int_0^t \e^{(t-s)A} Q_s(X_s)\d W_s$, which is an adapted continuous process on $\H$ up to the life time $\zeta$.  Then $Y_t:= X_t -\xi_t $ is the mild solution to   the equation
$$\d Y_t = \big(A Y_t + b_t (Y_t+\xi_t)\big)\d t,\ \ \ Y_0=X_0,\ \ t<\zeta.$$ Due to \eqref{C1} for $B=0$, the increasing property of $h,\Phi$,  and   $A\le 0$, this implies that for any $T>0$,
$$\d |Y_t|^2 \le 2\<Y_t,  b_t(Y_t+\xi_t)\>\d t\le 2\big(\Phi_{T\land \zeta}(|Y_t|^2)+ h_T(|\xi_t|)\big)\d t,\ \ |Y_0|^2=|X_0|^2, t<\zeta\land T.$$ Letting
\beq\label{AAA}\Psi_T(s)=\int_1^s\ff{\d r}{2\Phi_{T\land \zeta}(r)},\ \ \aa_T= |X_0|^2 + 2\int_0^{T\land\zeta} h_{T\land\zeta}(|\xi_s|)\d s, \ \ T> 0,\end{equation}  we obtain
$$|Y_t|^2 \le \aa_T +2\int_0^t \Phi_T(|Y_r|^2)\d r,\ \ \ T>0, t\in [0,T\land\zeta).$$By Biharis's inequality, this implies
\beq\label{B3} |Y_t|^2 \le \Psi_{T}^{-1} \big(\Psi_{T}(\aa_t)+t\big),\ \ T>0,   t\in [0,\zeta\land T).\end{equation} Moreover, {\bf (a1)} and   $\|Q\|_{T,\infty}<\infty$ yield
\beq\label{P*} \E\sup_{t\in [0,T\land\zeta)} |\xi_t|^2<\infty,\ \ T>0,\end{equation} so that on the set $ \{\zeta<\infty\}$ we have $\P$-a.s.
 \beq\label{P**} \limsup_{t\uparrow \zeta} |Y_t|= \limsup_{t\uparrow \zeta} |X_t|=\infty.\end{equation} We conclude that $\P(\zeta<\infty)=0$, i.e. $X_t$ is non-explosive. Indeed, on the set $\{\zeta\le T\}$,   \eqref{P*} implies $\P$-a.s.
 $$  \aa_T= |X_0|^2 +2 \int_0^\zeta h_\zeta (|\xi_s|)\d s <\infty,$$  so that  \eqref{B3}   and \eqref{P**}  imply
  $$\infty= \limsup_{t\uparrow\zeta} |Y_t|^2 \le   \Psi_{T}^{-1} \big(\Psi_{T}(\aa_T)+T\big)<\infty,$$ where the last step is due to the fact that $\Psi_T (r)\uparrow\infty$ as $r\uparrow\infty$, which implies  $\Psi_T^{-1}(r)<\infty$ for any $r\in (0,\infty).$  This    contradiction   means that  $\P(\zeta\le T)=0$ holds for all $T\in (0,\infty)$. Hence, $\P(\zeta <\infty)=0.$

  Finally, let $\aa_T(x)$ be defined in \eqref{AAA} as $\aa_T$ for $X_0=x$.  By the local boundedness of $\|Q_t\|_\infty,$   $\aa_T(x)$ is $\P$-a.s. locally bounded on $[0,\infty)\times \H$. Then, applying     \eqref{B3} to $X_0=y$ we conclude that      for any $x\in\H$ and $T>0$,   $\P$-a.s.
\beq\label{B4} \Xi^x:= \sup_{t\in [0,T], |y-x|\le 1}|X_t^y|<\infty,\end{equation}where $X_t^y$ is the mild solution for $X_0=y$.
Let $X_t^{(n,z)}$ solve \eqref{E1}   with $X_0=z$ for $B=0$ and  $(b^{[n]}, Q^{[n]})$ in place of $(b,Q)$, and let $$P_t^{(n)}f(z) =\E f(X_t^{(n,z)}),\ \ \ t\ge 0, f\in \B_b(\H), z\in\H, n\ge 1.$$ By Proposition \ref{P3.2}, $P_T^{(n)}$ is strong Feller. Since $X_t^{(n,z)}=X_t^z$ for $t\le \tau_n^z$, where $\tau_n^z:= n\land \inf\{t\ge 1: |X_t^{(n,z)}|\ge n\},$  it follows that
\beg{equation*}\beg{split} |P_T f(y)-P_Tf(x)| &\le |P_T^{(n)} f(y) -P_T^{(n)}f(x)| + 2 \|f\|_\infty \Big(\P(\tau_n^x\le T)+\P(\tau_n^y\le T)\Big) \\
&\le  |P_T^{(n)} f(y) -P_T^{(n)}f(x)| + 4\|f\|_\infty \P(\Xi^x\ge n),\ \ n>T, |y-x|\le 1.\end{split}\end{equation*} Since $P_T^{(n)}$ is strong Feller, this implies
$$\limsup_{y\to x} |P_T f(y)-P_Tf(x)|\le 4\|f\|_\infty \P(\Xi^x\ge n),\ \ \ n>T,\ \ f\in \B_b(\H).$$ Letting $n\to\infty$ and using \eqref{B4},  we obtain $\limsup_{y\to x} |P_T f(y)-P_Tf(x)|=0$ for $f\in \B_b(\H)$. Thus, $P_T$ is strong Feller.

\section{Proof of Theorem \ref{T1.2}}

Throughout this section, we assume {\bf (a1)},  {\bf (a2')} and {\bf (a3')}.
The idea of the proof is to transform \eqref{E1} into an equation with regular coefficients,  so that gradient estimates for the solution of the new equation can be derived. To this end, we use   the regularization representation \eqref{E4}. Let us fix    $T>0$. By Lemma \ref{L2.3}, we take large enough $\ll(T)>0$ such that for any $\ll\ge \ll(T)$ the unique solution  $u$ to \eqref{ABC} satisfies
\beq\label{O1}\|\nn^2 u\|_{T,\infty}+ \|\nn u\|_{T,\infty}\le \ff 1 8. \end{equation}   By \eqref{O1}, for any $t\in [0,T]$,
$$\H\ni x\mapsto \theta_t(x):= x+u_t(x)\in \H$$ is a  diffeomorphism with
\beq\label{O3} \ff 78 \le \|\nn \theta\|_{T,\infty}\le \ff 98,\ \ \ \ff 89\le \|\nn \theta^{-1}\|_{T,\infty}\le \ff 87.\end{equation}

Now, let $X_t^x$ solve \eqref{E1} for $X_0=x$. By \eqref{E4}, $Y_t^x:= \theta_t(X_t^x)$ satisfies
\beg{equation}\label{O2}\beg{split} Y_t^x= &\ \e^{tA}Y_0^x + \int_0^t \e^{(t-s)A} \big\{(\ll-A)u_s+B_s+\nn_{B_s}u_s\big\}\circ\theta_s^{-1}(Y_s^x)\d s \\
&+\int_0^t \e^{(t-s)A}\big\{Q_s+(\nn u_s)Q_s\big\}\circ\theta_s^{-1}(Y_s^x)\d W_s,\ \ t\in [0,T].\end{split}\end{equation} Thus, letting
$$\bar b_t= \{B_t+\nn_{B_t}u_t + (\ll-A) u_t\}\circ \theta_t^{-1},\ \ \bar Q_t= \{Q_t+(\nn u_t)Q_t\}\circ \theta_t^{-1},$$ $  \bar X_t^{x}:= Y_t^{\theta_0^{-1}(x)} $ is a  mild solution to the equation
\beq\label{O2'} \d \bar X_t^x= A\bar X_t^x\d t + \bar b_t(\bar X_t^x)\d t +\bar Q_t(\bar X_t^x)\d W_t,\ \ t\in [0,T], \bar X_0^x= x. \end{equation}
Let $\bar P_t f(x)= \E f(\bar X_t^x)$. We have
\beq\label{SM} \beg{split}  & P_t f(x):= \E f( X_t^x) =\E (f\circ \theta_t^{-1}) (Y_t^{x})\\
  &= \E (f\circ \theta_t^{-1}) (\bar X_t^{\theta_0(x)}) =   (\bar P_t f\circ\theta_t)(\theta_0(x)),\ \    f\in\B_b(\H), t\in [0,T], x\in\H.\end{split} \end{equation}  We first study gradient estimates and the log-Harnack inequality for $\bar P_t$. To this end,  one may wish to apply the corresponding results derived recently in \cite{WZ14}. However,  in the present case the assumption (A1) in \cite{WZ14} is not available, i.e. our conditions do not imply the existence of  $K\in L^2([0,T];\d t)$   such that
$$  |\e^{tA}(\bar b_s(x)-\bar b_s(y))| \le K(t) |x-y|,\ \ t,s\in [0,T], x,y\in\H.$$   Hence, we are not   at the position to apply results in \cite{WZ14}.

 To overcome the singularity of $\bar b$ caused by infinite-dimensions, we will use the finite-dimensional approximation argument. Unlike in the better situations of \cite{RW10, WZ14} where the desired gradient estimates and log-Harnack inequality have been established, in the present case we only have a weaker approximation result. More precisely,  letting $\bar X_t^{(n,x)}$ solving the finite-dimensional equation on $\H_n$:
\beq\label{FF} \d \bar X_t^{(n,x)} = \big\{ A\bar X_t^{(n,x)} + {\bar b}_t^{(n)}(\bar X_t^{(n,x)})\big\}\d t + {\bar Q}_t^{(n)}(\bar X_t^{(n,x)})\d W_t,\ \
\bar X_0^{(n,x)}=  x\in \H_n,t\in [0,T],\end{equation}
where $\bar b^{(n)}= \pi_n \bar b, \bar Q^{(n)}= \pi_n \bar Q$, instead of $\lim_{n\to 0} \E |\bar X_t^x-\bar X_t^{(n,\pi_n x)}|^2  =0$ for every $t\in [0,T],$ we only have
\beq\label{FF2} \lim_{n\to \infty} \int_0^T\E   |\bar X_t^x-\bar X_t^{(n,\pi_n x)}|^2\d t =0.\end{equation}
But this is already enough for our purpose.

\beg{lem}\label{L6.1} Assume {\bf (a1)}, {\bf (a2')} and {\bf (a3')}. For any $T>0$ and large enough $\ll\ge \ll(T)$, there exists a constant $C>0$   such that  the following assertions hold. \beg{enumerate}\item[$(1)$] If in addition $\|B\|_{T,\infty}<\infty$, then $\eqref{FF2}$ holds.
\item[$(2)$]  $|\nn \bar P_t^{(n)} f|^2 \le C \bar P_t^{(n)} |\nn f|^2,\ \ n\ge 1, t\in [0,T], f\in C_b^1(\H_n).$
\item[$(3)$]  $\ff t C |\nn \bar P_t^{(n)}f|^2 \le \bar P_t^{(n)} f^2 - (\bar P_t^{(n)} f)^2 \le Ct \bar P_t^{(n)} |\nn f|^2,\ \ n\ge 1, t\in [0,T], f\in C_b^1(\H_n). $
\item[$(4)$]  $\bar P_t^{(n)}\log f(y)\le \log \bar P_t^{(n)} f(x) +\ff {C|x-y|^2} t,  \ \ n\ge 1, t\in [0,T], 0<f\in \B_b(\H_n).$\end{enumerate}
\end{lem}

\beg{proof}  For simplicity, we omit $x$ and $\pi_n x$ from the subscripts, i.e. we write $(\bar X_t,\bar X_t^{(n)})$ instead of $(\bar X_t^x, \bar X_t^{(n,\pi_n x)}).$ The essential part of the proof is for (1) and (2), since   (3) and (4)  can be deduced  from (2) by using standard arguments.

 (1) By \eqref{O2'} and \eqref{FF}, we have
\beg{equation}\label{DC}\beg{split} &I_n := \E \int_0^T  \e^{-2p t} |\bar X_t-\bar X_t^{(n)}|^2\d t\\
&\le 5|x-\pi_n x|^2+  5\E\int_0^T \e^{-2p t} \bigg|\int_0^t \e^{(t-s)A} \big\{\bar b_s(\bar X_s)- \bar b_s^{(n)}(\bar X_s)\big\}\d s\bigg|^2\d t\\
&\quad +  5\E\int_0^T \e^{-2p t} \bigg|\int_0^t \e^{(t-s)A} \big\{\bar b_s^{(n)}(\bar X_s)- \bar b_s^{(n)}(\bar X_s^{(n)})\big\}\d s\bigg|^2\d t\\
&\quad +  5\E\int_0^T \e^{-2p t}  \int_0^t \big\|\e^{(t-s)A} \big\{\bar Q_s(\bar X_s)-\bar Q_s^{(n)}(\bar X_s)\big\}\big\|_{HS}^2 \d s  \\
&\quad +  5\E\int_0^T \e^{-2p t}  \int_0^t \big\|\e^{(t-s)A} \big\{\bar Q_s^{(n)}(\bar X_s)-\bar Q_s^{(n)}(\bar X_s^{(n)})\big\}\big\|_{HS}^2 \d s,\ \ p\ge 1.\end{split}\end{equation}
Obviously, {\bf (a1)}, {\bf (a2')} and {\bf (a3')} imply
$ \sup_{t\in [0,T],n\ge 1} \E\big(|X_t|^2+|\bar X_t^{(n)}|^2\big) <\infty, $ so that
\beq\label{NND} I:= \limsup_{n\to\infty} I_n <\infty.\end{equation}  Moreover, {\bf (a2')} and \eqref{O1} imply that  $\bar Q$ is bounded. So, it follows from    \eqref{O2'}, \eqref{O3}, $\sup_{t\in [0,T]}\E |X_t|^2<\infty$, and {\bf (a1)} that
\beg{equation*}\beg{split} & \sup_{t\in [0,T]}  \E\bigg|\int_0^t\e^{(t-s)A} \bar b_s (\bar X_s)\d s \bigg|^2\\
&\le 3\sup_{t\in [0,T]} \bigg\{\E |\bar X_t|^2+|\bar X_0|^2 +\E\bigg|\int_0^t \e^{(t-s)A}\bar Q_s(\bar X_s)\d W_s\bigg|^2\bigg\}\\
&\le 3\sup_{t\in [0,T]} \bigg\{\|\nn\theta\|_{T,\infty}^2  \E|X_t^{\theta_0^{-1}(x)}|^2 + |\theta_0^{-1}(x)|^2 +\|\bar Q\|_{T,\infty}^2   \int_0^t \|\e^{(t-s)A}\|_{HS}^2\d s \bigg\}<\infty.
\end{split}\end{equation*}  Thus, by the dominated convergence theorem,  \eqref{DC} implies
\beq\label{O4}\beg{split} I& \le  5 \limsup_{n\to\infty}\E\int_0^T \e^{-2p t} \bigg|\int_0^t \e^{(t-s)A} \big\{\bar b_s^{(n)}(\bar X_s)- \bar b_s^{(n)}(\bar X_s^{(n)})\big\}\d s\bigg|^2\d t\\
&+ 5\limsup_{n\to\infty}\E\int_0^T \e^{-2p t}  \int_0^t \big\|\e^{(t-s)A} \big\{\bar Q_s^{(n)}(\bar X_s)-\bar Q_s^{(n)}(\bar X_s^{(n)})\big\|_{HS}^2 \d s=: I'+I''.\end{split}\end{equation} By \eqref{00D} and {\bf (a2')} we have
\beq\label{O6}   \sup_{n\ge 1} \|\nn \bar Q^{(n)}\|_{T,\infty}^2 \le C_1 \end{equation} for some constant $C_1>0$,  so that
\beg{equation*}\beg{split} I'' &\le 5 C_1\limsup_{n\to\infty} \int_0^T \e^{-2p t} \d t \int_0^t \big\|\e^{(t-s)A}\big\|_{HS}^2 \E |\bar X_s- \bar X_s^{(n)}|^2 \d s \\
&\le 5 C_1  \limsup_{n\to\infty} \int_0^T \e^{-2p s}\E |\bar X_s-\bar X_s^{(n)}|^2\d s \int_s^T    \big\|\e^{(t-s)A}\big\|_{HS}^2 \e^{-2p(t-s)} \d t
 \le c(p) I,\end{split}\end{equation*} where according to {\bf (a1)},
 $$c(p):= 5C_1 \int_0^T \big\|\e^{tA}\big\|_{HS}^2 \e^{-2p t} \d t\to 0\ \ \text{as}\ p\to\infty.$$ Taking large enough $p>1$ such that $c(p)\le \ff 1 2,$ and substituting this into \eqref{O4}, we arrive at
 \beq\label{O9}\beg{split}  I&\le   10 \limsup_{n\to\infty} \sum_{i=1}^\infty \E\int_0^T \e^{-2p t} \bigg|\int_0^t \e^{-(t-s)\ll_i}\<\bar b_s(\bar X_s)-\bar b_s(\bar X_s^{(n)}), e_i\>\d s\bigg|^2\\
 &\le 10 \limsup_{n\to\infty} \sum_{i=1}^\infty \E\int_0^T \bigg(\int_0^t \e^{-(t-s)(p+\ll_i)}(p+\ll_i)\d s\bigg)\\
 &\qquad\qquad\qquad\times \bigg(\int_0^t\ff{\e^{-(t-s)(p+\ll_i)-2p s}}{p+\ll_i} \<\bar b_s(\bar X_s)-\bar b_s(\bar X_s^{(n)}), e_i\>^2\d s\bigg)\d t\\
 &\le 10 \limsup_{n\to\infty} \sum_{i=1}^\infty \E\int_0^T\bigg(\ff{\e^{-2p s}}{(p+\ll_i)^2} \<\bar b_s(\bar X_s)-\bar b_s(\bar X_s^{(n)}), e_i\>^2\bigg)\\
 &\qquad\qquad\qquad\times \bigg(\int_s^T \e^{-(t-s)(p+\ll_i)}(p+\ll_i)\d t\bigg)\d s\\
 &\le 10 \limsup_{n\to\infty} \sum_{i=1}^\infty \E\int_0^T\ff{\e^{-2p s}}{(p+\ll_i)^2} \<\bar b_s(\bar X_s)-\bar b_s(\bar X_s^{(n)}), e_i\>^2\d s.\end{split}\end{equation}
 Noting that $\<\bar b_s, e_i\>= (\ll+\ll_i)\<u_s\circ\theta_s^{-1},e_i\>+ \<(B_s+\nn_{B_s}u_s)\circ\theta_s^{-1}, e_i\>,$ for $p\ge \ll$ this implies
 \beq\label{O92} \beg{split}  &I \le 20 \limsup_{n\to\infty} \E\int_0^T \e^{-2p s} \big|u_s\circ\theta_s^{-1} (\bar X_s) - u_s\circ\theta_s^{-1} (\bar X_s^{(n)})\big|^2\d s\\
 &+ 20 \limsup_{n\to\infty} \E\int_0^T \e^{-2p s}\ff{|(B_s+\nn_{B_s}u_s)\circ\theta_s^{-1}(\bar X_s)- (B_s+\nn_{B_s}u_s)\circ\theta_s^{-1}(\bar X_s^{(n)})|^2}{p^2} \d s\\
 &\le 20 \Big(\big\|\nn(u\circ\theta^{-1})\big\|_{T,\infty}^2 +\ff{\|\nn\{(B+\nn_Bu)\circ\theta^{-1}\}\|_{T,\infty}^2}{p^2}\Big)I\le \ff 1 2I\end{split}
 \end{equation}
 for   large enough $p\ge \ll$, since due to \eqref{O1} and \eqref{O3} we have $\|\nn (u\circ\theta^{-1})\|_{T,\infty}\le \ff 1 7.$
 Combining this with \eqref{NND}, we prove $I=0$ which is equivalent to \eqref{FF2}.

 (2) By \eqref{FF},  \eqref{O6} and {\bf (a1)},    we have
 \beq\beg{split} \label{O7} \E| \bar X_t^{(n,x)}- \bar X_t^{(n,y)}|^2 \le &3 |x-y|^2 + 3 \bigg|\int_0^t \e^{(t-s)A} \big\{\bar b_s^{(n)}(\bar X_s^{(n,x)})-\bar b_s^{(n)}(\bar X_s^{(n,y)})\big\}\d s\bigg|^2\\
  &+C_2\int_0^t \E | \bar X_s^{(n,x)}- \bar X_s^{(n,y)}|^2 \d s,\ \ t\in [0,T], n\ge 1, x,y\in\H_n\end{split}\end{equation} for some constant $C_2>0$.  On the other hand, similarly to  \eqref{O9} and \eqref{O92}, for    large enough $p>0$,  there holds
 \beg{equation*}\beg{split} &3\E \int_0^T \e^{-2p t}\bigg|\int_0^t \e^{(t-s)A}\big\{\bar b_s^{(n)}(\bar X_s^{(n,x)})-\bar b_s^{(n)}(\bar X_s^{(n,y)})\big\} \d s\bigg|^2\d t\\
 &\le \ff 12 \int_0^T \e^{-2p t}\E| \bar X_t^{(n,x)}- \bar X_t^{(n,y)}|^2\d t,\ \ n\ge 1, x,y\in\H_n.\end{split}\end{equation*} Combining this with \eqref{O7}, we obtain
\beg{equation*}\beg{split} &\int_0^t \e^{-2p t}  \E| \bar X_t^{(n,x)}- \bar X_t^{(n,y)}|^2\d t\\
&\le 6\int_0^T \e^{-2p t} |x-y|^2 \d t + 2C_2 \int_0^T \e^{-2p t}\d t \int_0^t \E| \bar X_s^{(n,x)}- \bar X_s^{(n,y)}|^2 \d s\\
&\le 6 \int_0^T \e^{-2p t} |x-y|^2 \d t + 2C_2 \int_0^T \e^{-2p s} \E| \bar X_s^{(n,x)}- \bar X_s^{(n,y)}|^2 \d s\int_s^T \e^{-2p(t-s)}\d t \\
&\le  6\int_0^T \e^{-2p t} |x-y|^2 \d t +\ff {C_2}p \int_0^t \e^{-2p t} \E| \bar X_t^{(n,x)}- \bar X_t^{(n,y)}|^2\d t.\end{split}\end{equation*} Taking large enough $p_0=  2C_2$, such that $\ff{C_2}p\le \ff 1 2$ for $p\ge p_0$, we obtain
$$\int_0^T \e^{-2 p t} \E| \bar X_t^{(n,x)}- \bar X_t^{(n,y)}|^2\d t \le C |x-y|^2 \int_0^T \e^{-2pt}\d t,\ \ n\ge 1, x,y\in\H_n, p\ge p_0 $$   for some constant $C>0$. Since a finite measure on $[0,T]$ is determined by its Laplace transform, this implies that for   any $n\ge 1$ and  $x,y\in\H_n,$   $\E| \bar X_t^{(n,x)}- \bar X_t^{(n,y)}|^2\le C|x-y|^2$ holds for a.e. $t\in [0,T].$   By the continuity of $\E| \bar X_t^{(n,x)}- \bar X_t^{(n,y)}|^2$ in $t$, we prove
$$\E| \bar X_t^{(n,x)}- \bar X_t^{(n,y)}|^2\le C|x-y|^2,\ \ t\in [0,T], x,y\in\H_n, n\ge 1. $$   Then for any $f\in C_b^1(\H_n)$ and $(t,x)\in [0,T]\times \H_n$,
\beg{equation*}\beg{split} |\nn \bar P_t^{(n)} f(x)|^2&:=\limsup_{y\to x} \ff{|\bar P_t^{(n)}f(x)- \bar P_t^{(n)}f(y)|^2}{|x-y|^2} \le \limsup_{y\to x} \ff{\E|f(\bar X_t^{(n,x)}) -  f(\bar X_t^{(n,y)})|^2}{|x-y|^2}\\
&\le \big\{\bar P_t^{(n)} |\nn f|^2 (x)\big\} \limsup_{y\to x} \ff{\E |\bar X_t^{(n,x)}-X_t^{(n,y)}|^2}{|x-y|^2} \le C \bar P_t^{(n)} |\nn f|^2(x).\end{split}\end{equation*}

(3) By an approximation argument, we may and do assume that $\bar b^{(n)}\in C([0,T]; C_b^2(\H_n;\H_n))$ and $f\in C_0^2(\H_n)$. By {\bf (a2')} and \eqref{00D}, there exist    constants $c_1,c_2>0$ such that
\beq\label{ND} c_2I^{(n)}\ge  \bar Q_t^{(n)} (\bar Q_t^{(n)})^*\ge c_1 I^{(n)},\end{equation} where $I^{(n)}$ is the identity on $\H_n$. Let $\bar P_{s,t}^{(n)}$ be the in-homogenous Markov semigroup associated to \eqref{FF}.  We have
$$\bar P_t^{(n)}  f^2 -(\bar P_t^{(n)}f)^2 =\int_0^t\ff{\d}{\d s} \bar P_s^{(n)} (\bar P_{s,t}^{(n)} f)^2 \d s
= \int_0^t  \bar P_s^{(n)} \<\bar Q^{(n)}_s (\bar Q_s^{(n)})^*\nn  P_{s,t}^{(n)} f, \nn P_{s,t}^{(n)} f\> \d s.$$
Combining this with \eqref{ND} and (2), we prove (3). For instance, regarding $s$ as the starting time, we see that (2) also holds for $\bar P_{s,t+s}^{(n)}$ in place of $\bar P_t^{(n)}$, so that
\beg{equation*}\beg{split} &\bar P_t^{(n)}  f^2 -(\bar P_t^{(n)}f)^2 \le c_2 \int_0^t  \bar P_s^{(n)} |\nn  P_{s,t}^{(n)} f|^2   \d s\\
&\le c_2C \int_0^t \bar P_s^{(n)}\bar P_{s,t}^{(n)} |\nn f|^2\d s = c_2 C t \bar P_t^{(n)} |\nn f|^2.\end{split}\end{equation*}

(4) As in (3), we   assume that $\bar b^{(n)}\in C([0,T]; C_b^2(\H_n;\H_n)).$  It suffices to prove for $ f\in C_b^2(\H_n)$ which is strictly positive such that $\nn f=0$ outside a bounded set. Take $\ggm_s= x+ \ff s t (y-x), s\in [0,t].$ By \eqref{ND} and (2), we have
\beg{equation*} \beg{split} &\ff{\d}{\d s} \bar P_s^{(n)} \log \bar P_{s,t}^{(n)} f(\ggm_s)\\
& =-  P_s \<\bar Q^{(n)}_s (\bar Q_s^{(n)})^* \nn \log \bar P_{s,t}^{(n)} f, \nn \log \bar P_{s,t}^{(n)}f\> (\ggm_s) +\ff 1 t \<y-x, \nn \bar P_s^{(n)} \log \bar P_{s,t}^{(n)}\>  f(\ggm_s)\\
&\le \ff{|x-y|}t |\nn \bar P_s^{(n)} \log \bar P_{s,t}  f|(\ggm_s)- C_3\bar P_s^{(n)} |\nn \log \bar P_{s,t}^{(n)}f|^2(\ggm_s)\\
&\le \ff {|x-y|^2}{(2C_3t)^2},\ \ \ s\in [0,t] \end{split}\end{equation*} for some constant $C_3>0$. Integrating over $[0,t]$ we prove (4) for some constant $C>0$.
\end{proof}

\beg{proof}[Proof of Theorem \ref{T1.2}] (a) We first assume that $\|B\|_{t,\infty}<\infty$ for $T>0.$   In this case we   observe  that due to Lemma \ref{L6.1}(1), assertions in Lemma \ref{L6.1} (2)-(4) hold for $\bar P_t$ in place of $\bar P_t^{(n)}.$
To save space we only prove the first inequality in (3), the proofs for others are similar and even simpler. Let $f\in C_b^1(\H).$ Since $\bar P_t 1=1,$ this inequality is equivalent to
$|\nn\bar P_t f|^2\le \ff C t P_t f^2.$  Since $\bar P_t f\in C_b(\H)$ which  is true even for $f\in \B_b(\H)$ according to the strong Feller property of $P_t$ and the relation \eqref{SM},     this inequality follows from
\beq\label{LST} \ff{|\bar P_t f(x) -\bar P_t f(y)|^2}{|x-y|^2}\le \ff C t \int_0^1 \bar P_t f^2 (x+s(y-x))\d s.\end{equation}
By Lemma \ref{L6.1}(3), we have
\beg{equation}\label{LD} \beg{split}  \ff{|(\bar P_t^{(n)} f)(\pi_n x) -(\bar P_t^{(n)} f)(\pi_n y)|^2}{|x-y|^2}&\le \bigg(\int_0^1 |\nn \bar P_t^{(n)} f|\circ\pi_n (x+s(y-x)) \d s\bigg)^2\\
&\le  \ff C t \int_0^1( \bar P_t^{(n)} f^2)\circ\pi_n (x+s(y-x))\d s. \end{split}\end{equation}Moreover, Lemma \ref{L6.1}(1) implies
$$\lim_{n\to\infty}  \int_0^T |\bar P_t f - (\bar P_t^{(n)} f)\circ\pi_n|^2\d t =0,\ \ f\in C_b^1(\H).$$ So, by letting $n\to\infty$ (up to a subsequence) in \eqref{LD} we prove \eqref{LST} for a.e.   $t\in [0,T]$ with respect to the Lebesgue measure. By the continuity of $\bar X_t$, $\bar P_t f$ is continuous in $t$. Therefore, \eqref{LST} holds for all $t\in [0,T].$

Now, according to \eqref{SM} and \eqref{O3}, we only need to prove Theorem \ref{T1.2} for $\bar P_t f$ in place of $P_t$.  By the above observation,
Theorem \ref{T1.2}(1) as well as \eqref{G0} and \eqref{LH0} with $t\in (0,1]$ hold for $\bar P_t$ in place of $P_t$.
 Then the proof is complete by the following two facts:   (a) Due to the semigroup property and Jensen's inequality, if \eqref{G0} and \eqref{LH0} hold for $t\in (0,1]$, then they also hold for   all $t>0$; (b) According to \cite[Proposition 1.3]{W14}, \eqref{G0} is equivalent to \eqref{H0}.

(b) In general, let $\tt P_t^{(n)}$ be the semigroup associated to \eqref{E1} for $\tt B^{(n)}:= B\circ\psi_n$ in place of $B$, where
$$\psi_n(x):= \Big(1\land \ff n{|x|}\Big) x,\ \ x\in\H.$$ Then $\|\tt B^{(n)}\|_{T,\infty}<\infty$ for $T>0$ and   {\bf (a2')} holds for $\tt B^{(n)}$ in place of $B$ with the same function $\Psi$. According to (a), assertions in Lemma \ref{L6.1} (2)-(4) hold for $\tt P_t^{(n)}$ in place of $\bar P_t^{(n)}.$
Moreover, by the uniqueness and non-explosion of solutions to the equation \eqref{E1}, we have
$$\lim_{n\to\infty} \tt P_t^{(n)} f= P_t f,\ \ t\ge 0, f\in C_b(\H).$$  Therefore, as explained in (a) that these assertions also hold for $P_t$.   \end{proof}

\paragraph{\bf Acknowledgement.} The author would like to thank  Jian Wang,  X. Huang  and   X. Zhang  for helpful comments and corrections.

\beg{thebibliography}{99}

\bibitem{AG} A. Alabert, I. Gy\"ongy, \emph{On stochastic reaction-diffusion equations with singular force term,} Bernoulli 7(2001), 145--164.

\bibitem{ATW14} M. Arnaudon, A. Thalmaier, F.-Y. Wang, \emph{Equivalent Harnack and gradient inequalities for pointwise curvature lower bound,} to appear in Bull. Sci. Math. DOI: 10.1016/j.bulsci.2013.11.001.

\bibitem{Bakry} D. Bakry, I. Gentil, M. Ledoux, \emph{Analysis and Geometry of Markov Diffusion Operators,} Springer, Berlin, 2014.

\bibitem{BGL} D. Bakry, I. Gentil, M. Ledoux, \emph{On Harnack inequalities and optimal transportation,} arXiv: 1210.4650. To appear in Annali della Scuola Normale Superiore di Pisa.


\bibitem{DF}  G. Da Prato, F. Flandoli, \emph{Pathwise uniqueness for  a class of SDE in Hilbert spaces and applications,} J. Funct. Anal. 259(2010), 243--267.


\bibitem{DR1}   G. Da Prato, F. Flandoli, E. Priola, M. R\"ockner, \emph{Strong uniqueness for stochastic evolution equations in Hilbert spaces perturbed by a bounded measurable drift,}
  Ann. Probab. 41(2013), 3306--3344.

\bibitem{DR2}  G. Da Prato, F. Flandoli, E. Priola, M. R\"ockner, \emph{Strong uniqueness for stochastic evolution equations with  unbounded measurable drift term,} to appear in J. Theor. Probab. DOI: 10.1007/s10959-014-0545-0.

\bibitem{RN}   G. Da Prato, F. Flandoli, M. R\"ockner,  A. Yu. Veretennikov, \emph{Strong uniqueness for SDEs in Hilbert spaces with non-regular drift,}
  arXiv: 1404.5418.

\bibitem{DZ}  G. Da Prato, J. Zabczyk, \emph{Stochastic Equations in Infinite Dimensions,}  Cambridge University Press, Cambridge, 1992.


\bibitem{FP} F. Flandoli, M. Gubinelli, E. Priola, \emph{Well-posedness of transport equation by stochastic perturbation,} Invent. Math. 180(2010), 1--53.

\bibitem{FGP}  F. Flandoli, M. Gubinelli, E. Priola, \emph{Flow of diffeomorphisms for SDEs with unbounded H\"older continuous drift,} Bull. Sci. Math. 134(2010), 405--422.

\bibitem{GP} I. Gy\"ongy, E. Pardoux, \emph{On the regularization effect of space-time white noise on quasi-linear parabolic partial differential equations,} Probab. Theory Relat. Fields 97(1993), 211--229.

\bibitem{KR} N. V. Krylov, M. R\"ockner, \emph{Strong solutions of stochastic equations with singular time dependent drift,} 131(2005), 154--196.

\bibitem{LLW} H. Li, D. Luo, J. Wang, \emph{Harnack inequalities for SDEs with H\"older continuous drift,}  arXiv:1310.4382.

\bibitem{OD04} M. Ondrej\'et, \emph{Uniqueness for stochastic evolution equations in Banach spaces,}  Dissertationes Math. (Rozprawy Mat.) 426(2004).

\bibitem{LR13}  C. Pr\'ev\^ot,  M. R\"ockner, \emph{A Concise Course on Stochastic   Partial Differential Equations,} Lecture Notes in Math. Vol. 1905, Springer, Berlin, 2007.

\bibitem{RW10} 	M. R\"ockner, F.-Y. Wang, Log-Harnack  Inequality for Stochastic differential equations in Hilbert spaces and its consequences, Infinite Dimensional Analysis, Quantum Probability and Related Topics 13(2010), 27--37.

\bibitem{W97}   F.-Y. Wang, \emph{Logarithmic Sobolev
inequalities on noncompact Riemannian manifolds,} Probab.
Theory Related Fields 109(1997), 417--424.

\bibitem{W10} F.-Y. Wang, \emph{Harnack inequalities on manifolds with boundary and applications,}    J.
Math. Pures Appl.     94(2010), 304--321.

\bibitem{Wbook} F.-Y. Wang,  \emph{Harnack Inequalities for Stochastic Partial Differential Equations,} Springer, Berlin, 2013.

\bibitem{Wb2} F.-Y. Wang,  \emph{Analysis for Diffusion Processes on Riemannian Manifolds,} World Scientific, Singapore, 2013.

\bibitem{W14}   F.-Y. Wang, \emph{Derivative formula and gradient estimates for  Gruschin type semigroups,} to appear in J. Theo. Probab. DOI: 10.1007/s10959-012-0427-2.

\bibitem{WZ14}  	F.-Y. Wang, T. Zhang, \emph{Log-Harnack inequalities for semi-linear SPDE with strongly multiplicative noise,} Stoch. Proc. Appl. 124(2014), 1261--1274.

\bibitem{YW} T. Yamada, S. Watanabe, \emph{On the uniqueness of solutions of stochastic differential equations,} J. Math. Kyoto Univ. 11(1971), 155--167.

\bibitem{V} A. J. Veretennikov, \emph{Strong solutions and explicit formulas for solutions of stochastic integral equations,} (Russian) Mat. Sb. (N.S.) 111(480)(1980), 434--452.

\bibitem{Zh} X. Zhang, \emph{Strong solutions of SDES with singular drift and Sobolev diffusion coefficients,} Stoch. Proc. Appl. 115(2005), 1805--1818.

\bibitem{Z} A. K. Zvonkin, \emph{A transformation of the phase space of a diffusion process that will remove the drift,} (Russian) Mat. Sb. (N.S.) 93(152)(1974), 129--149.
\end{thebibliography}
\end{document}